# Homotopy spectral sequences (*)

Marco Grandis

Dipartimento di Matematica, Università di Genova, Via Dodecaneso 35, 16146-Genova, Italy

**Abstract.** In homotopy theory, exact sequences and spectral sequences consist of groups and pointed sets, linked by actions. We prove that the theory of such exact and spectral sequences can be established in a categorical setting which is based on the existence of kernels and cokernels *with respect to an assigned ideal of null morphisms*, a generalisation of abelian categories and Puppe-exact categories.

## 0. Introduction

The purpose of this paper is to study the homological aspects of exact sequences and spectral sequences coming from (unstable) homotopy theory; for instance the homotopy spectral sequence of a tower of fibrations (Bousfield and Kan [BK], IX.4) or of cofibrations (Baues [Bau], III.2). Such sequences of groups degenerate, in low dimension, into pointed sets, possibly acted on by groups.

To take this into account, we will use a categorical setting, the notion of 'homological category', which was introduced by the author at the Conference in Category Theory, Como 1990 and studied in the early 1990's [G6-G9]. It is a generalisation of Puppe-exact categories [Pu, Mt, HS], based on the existence of kernels and cokernels *with respect to an assigned ideal of null morphisms* and on a 'homology axiom'; the latter allows one to deal with subquotients (like homology, or the terms of a spectral sequence) and the induced morphisms between them.

Notice that the existence of categorical products is not assumed – a fact which is crucial for studying the coherence of induced morphisms, as showed in three previous papers on 'distributive homological algebra' [G2-G5]. The main result of this approach was the construction of universal models for spectral sequences *in Puppe-exact categories*. While this is (obviously) sufficient for the abelian context, the present extension will hopefully set the bases for universal models in *homological categories*, in order to cover the spectral sequences coming from homotopy theory.

Therefore, the hierarchy of categorical settings used here (and recalled in Section 1):

- *semiexact, homological, generalised exact, Puppe-exact category*,

should not be confused with a different, well-known system based on the existence of finite limits:

---

(*) Work supported by grants of Università di Genova and INDAM (Italy).

- *protomodular* (Bourn [B1]), *Borceux-Bourn homological* [BB, BC], *Barr-exact* [Ba], *semiabelian category* (in the sense of Janelidze, Márki and Tholen [JMT]).

(There is also a notion of 'homological monoid', by Hilton and Ledermann [HL], which is a pointed category satisfying three self-dual axioms; in particular, the category of groups is of this kind. The author is indebted to G. Janelidze for this reference.)

This second approach is appropriate to investigate other non-abelian properties, where the existence of finite limits is important. An interesting recent result in this framework is the 'semi-abelian Dold-Kan theorem' of D. Bourn [B2]: if **C** is a semi-abelian category, then the Moore normalisation functor M: Simpl**C** → Ch**C** is monadic; this allows one to recover simplicial objects on **C** as Eilenberg-Moore algebras over chain complexes on **C**.

It would be good to have a clearer understanding of the cleavage between these two approaches, which is somewhat related to the cleavage between projective and affine geometry. As briefly discussed in 1.8, the 'perspective', or 'pre-projective' character of the present approach appears:

- in the crucial use of lattices of (normal) subobjects and their direct or inverse images,

- in the construction of the 'projective category' associated to a Puppe-exact category [G1], in its extension to homological categories given here (1.8), and in the characterisation of the projective categories associated to the abelian ones, established in [CG].

This paper was exposed at the 'International Conference on Category Theory', Cape Town, 2009. The author acknowledges various useful discussions with George Janelidze, Zurab Janelidze and Dominique Bourn.

*Outline*. In the first two sections we review the notions of semiexact and homological category, introduced and studied in [G6-G9]. A basic example is the category **Set**$_2$ of pairs of sets (2.1), which is homological with respect to the ideal of maps f: $(X, X_0) \to (Y, Y_0)$ such that $f(X) \subset Y_0$; it is a monoidal closed category over which every semiexact category is naturally enriched. The category **Top**$_2$ of pairs of topological spaces has a similar homological structure; homology theories in the sense of Eilenberg-Steenrod form a connected sequence of functors on it. Similarly, the category **Gp**$_2$ of pairs of groups is homological, and the natural framework of relative cohomology of groups.

In Section 3, the Massey exact couple and its associated spectral sequence are extended to homological categories. The next section shows that the homotopy spectral sequence of a tower of fibrations *between path-connected spaces* can be seen in the present framework of homological categories, as produced by an exact couple living in the homological category **Ngp** of *normalised groups* (4.1), a category of fractions of **Gp**$_2$ in which **Gp** embeds, in such a way that all the homomorphisms of groups become exact morphisms of **Ngp**. This category **Ngp** has a simple description as a quotient of a concrete homological category (4.5), inspired by a category of *homogeneous spaces* introduced by Lavendhomme [La].

In Sections 5 and 6 we remove the restriction to pathwise connected spaces, using the category **Act** of actions of groups on pointed sets, or actions for short, in which **Set.**, **Gp** and **Gp**$_2$ embed naturally. We prove that **Act** is homological and show that the homotopy spectral sequence of a tower of fibrations can be established in a suitable homological category **Nac** of *normalised actions*; this is a category of fractions of **Act**, which again can be realised as a quotient of a concrete homological category.

A universe $\mathcal{U}$ is fixed throughout; its elements are called small sets.

## 1. Semiexact and homological categories

We review, without proofs, the notions of semiexact, homological and generalised exact category, introduced in [G6-G9] as an extension of Puppe-exact categories [Pu, Mt, HS]. Notice that a condition of 'normal well-poweredness' of such categories is assumed, from Section 1.5 on.

**1.1. Semiexact categories.** A *semiexact* category is a pair $(\mathbf{A}, \mathcal{N})$, generally denoted as $\mathbf{A}$, satisfying the following two axioms (explained below):

(ex0)  $\mathbf{A}$ is a category and $\mathcal{N}$ is a *closed ideal* of $\mathbf{A}$,

(ex1)  every morphism $f: A \to B$ of $\mathbf{A}$ has a kernel and cokernel, with respect to $\mathcal{N}$.

The morphisms of $\mathcal{N}$ are called *null morphisms* of $\mathbf{A}$; the objects whose identity is null are called *null objects* of $\mathbf{A}$. (The condition that $\mathcal{N}$ be an ideal means that every composition with a null morphism yields a null one, while the closedness of $\mathcal{N}$ means that every null morphism factorises through a null identity. Equivalently, we can assign a set of null objects, closed under retracts.)

The kernel and cokernel of $f: A \to B$ will be written as:

(1)   ker f: Ker f $\rightarrowtail$ A,                    cok f: B $\twoheadrightarrow$ Cok f.

The kernel is characterised up to isomorphism by a universal property, *with respect to* $\mathcal{N}$:

(2)   the composite f∘ker f  is null,

    for every map  a  in  $\mathbf{A}$,  if  fa  is (legitimate and) null, then  a  factorises uniquely through  ker f.

The cokernel is defined by the dual property. Therefore the morphism ker f is mono, and cok f is epi. A *normal mono* is a kernel (of some morphism), and a *normal epi* is a cokernel; the arrows $\rightarrowtail$ , $\twoheadrightarrow$  *are reserved for such morphisms*.

Kernels and cokernels establish an anti-isomorphism between the ordered classes of normal subobjects and normal quotients of an object A.  They also yield the *normal coimage* and the *normal image* of the morphism  f: A $\to$ B:

(3)   ncm f  =  cok ker f,                    nim f  =  ker cok f,

which are respectively the least normal quotient of  A  and the least normal subobject of  B  through which  f  factorises.

A semiexact category is said to be *pointed semiexact,* or *p-semiexact,* if it has a zero object *and* its null objects coincide with the zero objects (all isomorphic, of course).

**1.2. Normal factorisations.** In a semiexact category $\mathbf{A}$, we say that the morphism f is N-mono if ker f is null, or equivalently if: $fh \in \mathcal{N}$ implies $h \in \mathcal{N}$; dually, f is N-epi if cok f is null. One can prove ([G6], 3.7) that every normal mono is N-mono and every normal epi is N-epi; this property is not trivial, and depends on the closedness of the ideal $\mathcal{N}$.

(As a counterexample, consider the selfdual category associated to the ordered set $\mathbb{Z}$ of integers, equipped with the *non-closed* ideal of strict inequalities $x < y$. Kernels exist: $\ker 1_x = (x{-}1 < x)$ and $\ker(x < y) = 1_x$; analogously for cokernels; but $\ker(\ker 1_x) = 1_{x-1}$ is not null.)

It follows easily that a morphism $f: A \to B$ factorises through its normal coimage and its normal image, by a unique morphism $g$

(1)
$$\begin{array}{ccccccc} \text{Ker } f & \rightarrowtail & A & \xrightarrow{f} & B & \twoheadrightarrow & \text{Cok } f \\ & & p \downarrow & & \uparrow m & & \\ & & \text{Ncm } f & \xrightarrow{g} & \text{Nim } f & & \end{array}$$

This is called *the normal factorisation*, $f = mgp$. We say that $f$ is an *exact* morphism if $g$ is an isomorphism; this happens if and only if $f$ can be factorised as a normal epi followed by a normal mono.

A semiexact category **A** is called *generalised exact* if all its morphisms are exact. An exact category which is p-semiexact (1.1) is the same as a Puppe-exact category, and will also be called a *p-exact* category.

**1.3. Exact sequences and exact functors.** In a semiexact category **A**, the sequence (f, g)

(1) $\qquad A \xrightarrow{f} B \xrightarrow{g} C \qquad\qquad A \xrightarrow{f} B \xrightarrow{g} C$
$\qquad\qquad\qquad\qquad\qquad\qquad\qquad\quad {}^{f'}\nearrow \quad\; {}^{g'}\searrow$
$\qquad\qquad\qquad\qquad\qquad\qquad\qquad\; \bullet \qquad\qquad \bullet$

is said to be *of order two* if $gf$ is null (iff $\text{nim } f \leq \ker g$, iff $\text{cok } f \geq \text{ncm } g$). It is said to be *exact* (in the object B) if it satisfies the following equivalent conditions:

(i) $\text{nim } f = \ker g$,

(ii) $\text{cok } f = \text{ncm } g$,

(iii) for every f', g' as in the right diagram above, if gf' and g'f are null, so is g'f'.

More particularly, the sequence (f, g) is said to be *short exact* if $f = \ker g$ and $g = \text{cok } f$.

A functor F: **A** → **B** between semiexact categories will be said to be *exact* if it preserves kernels and cokernels. Then F preserves null morphisms (f is null if and only if ker f = 1), null objects, normal factorisations, exact morphisms, exact and short exact sequences.

More generally, we say that F: **A** → **B** is:

- an *N-functor* if it preserves null morphisms,

- *left exact* if it preserves kernels,

- *right exact* if it preserves cokernels,

- *short exact* if it preserves short exact sequences,

- *long exact* if it preserves exact sequences.

It is easy to see that a functor is exact if and only if it is both short exact *and* long exact. Moreover, if F: **A** $\rightleftarrows$ **B** :G is an adjunction between semiexact categories *and* F, G are N-functors, then the left adjoint F is right exact, while the right adjoint G is left exact.

A *semiexact subcategory* of the semiexact category **A** is a subcategory **A**' satisfying:

(a) for every f in **A**' there is *some* kernel and *some* cokernel of f in **A** which belong to **A**',

(b) if m is a normal mono of **A** which belongs to **A**', and mf is in **A**', so is f; dually, if p is a normal epi of **A** which belongs to **A**', and fp is in **A**', so is f.

Then **A**', equipped with the ideal $\mathcal{N}' = \mathbf{A}' \cap \mathcal{N}$, is a semiexact category; further, its inclusion in **A** is exact and *conservative* (reflects the isomorphisms).

**1.4. Lattices and connections.** The category **Ltc** of *lattices and connections* formalises the structure of normal subobjects in semiexact categories, together with their direct and inverse images.

An object is a small lattice (always assumed to have 0 and 1). A morphism $f = (f_\bullet, f^\bullet): X \to Y$, called a *connection*, is a 'covariant Galois connection' between the lattices X and Y, i.e. an adjunction $f_\bullet \dashv f^\bullet$:

(i) $f_\bullet: X \to Y$ and $f^\bullet: Y \to X$ are increasing mappings,

(ii) $f^\bullet f_\bullet \geq \text{id } X$, $f_\bullet f^\bullet \leq \text{id } Y$ (in the category of ordered sets and increasing mappings).

Therefore, $f_\bullet$ preserves the existing joins (including 0), $f^\bullet$ preserves the existing meets (including 1) and:

(1) $f_\bullet f^\bullet f_\bullet = f_\bullet$, $\quad\quad\quad\quad\quad\quad\quad\quad f^\bullet f_\bullet f^\bullet = f^\bullet$,

(2) $f^\bullet(y) = \max\{x \in X \mid f_\bullet(x) \leq y\}$, $\quad\quad f_\bullet(x) = \min\{y \in Y \mid f^\bullet(y) \geq x\}$.

The composition is obvious. Isomorphisms can be identified with ordinary lattice-isomorphisms. The category **Ltc** is *selfdual*, under the contravariant endofunctor which carries each lattice to the opposite one and reverses any connection.

**Ltc** has a zero object: the one-point lattice $\mathbf{0} = \{*\}$, and the zero-morphism $f: X \to Y$ is described by: $f_\bullet(x) = 0$, $f^\bullet(y) = 1$. Kernels and cokernels exist, and the normal factorisation of f is:

(3)
$$\begin{array}{ccccccc} & & m & & f & & p \\ \downarrow f^\bullet 0 & \rightarrowtail & X & \longrightarrow & Y & \twoheadrightarrow & \uparrow f_\bullet 1 \\ & & q \downarrow & & \uparrow n & & \\ & & \uparrow f^\bullet 0 & \xrightarrow{g} & \downarrow f_\bullet 1 & & \end{array}$$

$m_\bullet(x) = x$, $\quad\quad m^\bullet(x) = x \wedge f^\bullet 0$, $\quad\quad p_\bullet(y) = y \vee f_\bullet 1$, $\quad p^\bullet(y) = y$,

$q_\bullet(x) = x \vee f^\bullet 0$, $\quad q^\bullet(x) = x$, $\quad\quad\quad\quad n_\bullet(y) = y$, $\quad\quad\quad\quad n^\bullet(y) = y \wedge f_\bullet 1$,

$g_\bullet(x) = f_\bullet(x)$, $\quad\quad g^\bullet(y) = f^\bullet(y)$.

Therefore, the morphism f is exact if and only if, for every $x \in X$ and $y \in Y$:

(a) $f^\bullet f_\bullet(x) = x \vee f^\bullet 0$, $\quad\quad\quad\quad\quad\quad f_\bullet f^\bullet(y) = y \wedge f_\bullet 1$.

Every element $a \in X$ determines a short exact sequence in **Ltc**:

$$\text{(4)} \quad \downarrow a \overset{m}{\rightarrowtail} X \overset{p}{\twoheadrightarrow} \uparrow a, \qquad (m_\bullet(x) = x, \quad m^\bullet(x) = x \wedge a; \qquad p_\bullet(x) = x \vee a, \quad p^\bullet(x) = x),$$

and all short exact sequences are of this type, up to isomorphism.

We are also interested in the (p-exact) subcategory **Mlc** of *modular lattices and modular connections* [G1], where a morphism is any exact connection between modular lattices. The latter are characterised by the conditions (a) above; but the modularity of the lattices makes these conditions equivalent to the following ones, which are plainly stable under composition:

(b) $\quad f^\bullet(f_\bullet x \vee y) = x \vee f^\bullet y, \qquad f_\bullet(f^\bullet y \wedge x) = y \wedge f_\bullet x \qquad \text{(for } x \in X, y \in Y\text{)}.$

It is interesting to note that the category **Ltc** is *semiadditive*, i.e. it has finite biproducts; this defines in the usual way a sum of parallel maps, which lacks opposites.

Indeed, if X and Y are lattices, their cartesian product X×Y (i.e., their product in the category of homomorphisms of lattices) is their *biproduct* in **Ltc**, by means of the following connections where the injections i, j are normal monos and the projections p, q are normal epis

$$\text{(5)} \quad X \; \underset{p}{\overset{i}{\rightleftarrows}} \; X \times Y \; \underset{q}{\overset{j}{\leftrightarrows}} \; Y \qquad\qquad i_\bullet \dashv i^\bullet = p_\bullet \dashv p^\bullet,$$

$$i_\bullet(x) = (x, 0_Y), \qquad i^\bullet(x, y) = x = p_\bullet(x, y), \qquad p^\bullet(x) = (x, 1_Y).$$

Now, the sum f+g of parallel maps X → Y is defined by means of the diagonal of X and of the codiagonal of Y

$$\text{(6)} \quad X \longrightarrow X \times X \overset{f \times g}{\longrightarrow} Y \times Y \longrightarrow Y \qquad\qquad f + g = (f_\bullet \vee g_\bullet, f^\bullet \wedge g^\bullet).$$

This sum is idempotent, f+f = f, and our structure is not additive. Actually, f+g is the join f∨g with respect to an obvious order relation between parallel connections, given by $f_\bullet \leq g_\bullet$ (or equivalently $f^\bullet \geq g^\bullet$).

The subcategory **Mlc** does not inherit the semiadditive structure of **Ltc**: if X and Y are modular lattices, the maps of the biproduct-decomposition of X×Y, in (5), are modular connections but the universal property for the product or sum is not satisfied.

**1.5. Perspective functor and modularity.** In the semiexact category **A**, each object A has a (possibly large) lattice NsbA of normal subobjects and a lattice NqtA of normal quotients, anti-isomorphic via kernels and cokernels.

Each morphism f: A → B has *direct and inverse images* for normal subobjects:

(1) $\quad f_*: \text{Nsb}A \to \text{Nsb}B, \qquad\qquad f_*(x) = \text{nim } fx = \ker \text{cok}(fx),$

$\quad\;\; f^*: \text{Nsb}B \to \text{Nsb}A, \qquad\qquad f^*(y) = \ker((\text{cok } y) \circ f) = \text{pullback of } y \text{ along } f,$

which form a Galois connection $f_* \dashv f^*$. In particular, if m: M ↣ A and p: A ↠ P

(2) $\quad m_*(m^*(x)) = x \wedge m, \qquad\qquad p^*(p_*(x)) = x \vee \ker(p) \qquad\qquad (x \in \text{Nsb}A).$

Analogously there are direct and inverse images for normal quotients, linked to the connection (1) by the ker-cok anti-isomorphism.

We say that **A** is a *semiexact U-category* if it is semiexact, all its objects and morphisms belong to U and moreover all its lattices NsbA of normal subobjects belong to U. As a consequence, all the anti-isomorphic lattices of normal quotients are also small. (Notice that the set of objects, the set of morphisms and even the *hom-sets* of **A** are *not* assumed to belong to U.) **Ltc** and **Mlc** satisfy our condition (cf. 1.4.4).

Since we shall mainly consider such categories, *semiexact category will mean, from now on, semiexact U-category*, while *unrestricted* semiexact category will refer to the original defintition 1.1.

Therefore, a semiexact category has a *perspective functor*, or *transfer functor for normal subobjects*, with values in the category of small lattices and connections (1.4):

(3)  $\text{Nsb}_\mathbf{A}: \mathbf{A} \to \mathbf{Ltc}$, $\qquad\qquad$ $A \mapsto \text{Nsb}A$, $\quad$ $f \mapsto (f_*, f^*)$,

which generally is just *long exact*.

One can prove (cf. [G6], 5.7) that the perspective functor is exact if and only if **A** is an *ex2-category*, i.e. it also satisfies the axiom:

(ex2)  normal monomorphisms and normal epimorphisms are stable under composition.

An object A of the semiexact category **A** is said to be *modular* if its lattice NsbA of its normal subobjects is modular. A morphism f: A → B is *left modular*, or *right modular*, or *modular* if its associated connection Nsb(f): NsbA → NsbB satisfies (4), or (5), or both:

(4)  $f^* f_* x = x \vee f^* 0$ $\qquad$ (for $x \in \text{Nsb}A$),
(5)  $f_* f^* y = y \wedge f_* 1$ $\qquad$ (for $y \in \text{Nsb}B$).

More particularly, we say that f is *left modular on* x (resp. *right modular on* y) when the above property holds for this particular normal subobject x of A (resp. y of B).

A semiexact category **A** is said to be *modular* if all its objects and morphisms are modular, or equivalently if its perspective functor Nsb: **A** → **Ltc** takes values in the (p-exact) subcategory **Mlc** of modular lattices and modular connections (1.4). One can prove that a modular semiexact category is also ex2, and that every generalised exact category is modular.

An exact functor F: **A** → **B** between semiexact categories defines a family of lattice *homomorphisms* (for A in **A**):

(6)  $(\text{Nsb } F)_A: \text{Nsb}_\mathbf{A}(A) \to \text{Nsb}_\mathbf{B}(FA)$, $\qquad\qquad$ $x \mapsto \text{nim}(Fx)$.

(Notice that Fx is a normal mono of codomain FA, whereas nim(Fx) is *the* normal *subobject* equivalent to the former.)

We say that F is *nsb-faithful* if all these homomorphisms are injective: if x, x' are normal subobjects of any A and Fx ~ Fx', then x = x'. We say that F is *nsb-full* if these homomorphisms are surjective: if y is a normal subobject of FA there exists some normal subobject x of A such that Fx ~ y.

**1.6. Homological categories.** A *homological* category will be an ex2-category **A** (1.5) satisfying the following *subquotient* axiom, or *homology* axiom:

(ex3) given a normal mono m: M ↣ A and a normal epi q: A ↠ Q, with m ≥ ker q (cok m ≤ q), the morphism qm is exact.

Then the normal factorisation of qm produces a commutative square, which can be proved to be *bicartesian*, i.e. pullback and pushout ([G6], 6.2):

(1)
$$\begin{array}{ccc} M & \xrightarrow{m} & A \\ h \downarrow & & \downarrow q \\ S & \xrightarrow{k} & Q \end{array}$$
(m ≥ ker q).

The object S ≅ Ncm(qm) ≅ Nim(qm) will be said to be a *subquotient* of A, and written as M/N, where N = Ker q. *The notation* s: M/N ↠ A *will refer to this bicartesian square.* (If **A** is p-exact, s can be viewed as a monorelation from M/N to A, i.e. a monomorphism in the category of relations on **A**.) The *numerator* and *denominator* of the subquotient S = M/N are the following subobjects of A

(2)    Num(S) = M,                Den(S) = Ker q ≤ Num(S),

under the usual abuse of notation: writing the object M (or Ker q) for the subobject m (or ker q), and S for the subquotient s.

A semiexact subcategory (1.3) of a homological category is homological ([G6], 6.4). The category **Ltc** is pointed homological, not exact; its subcategory **Mlc** is p-exact.

This notion of homological category, intermediate between semiexact and generalised exact category, has been studied in [G9, G10] and, in a more complete form, in [G6, G7]. It allows one to develop a consistent part of homological algebra, including the study of spectral sequences.

**1.7. Regular induction.** Let us have a morphism f: A → B in a homological category **A** and two subquotients s: M/N ↠ A, t: H/K ↠ B. We say that f *has a regular induction* from M/N to H/K if (with, again, the usual abuse of notation for subobjects)

(1)    $f_*(M) \le H$,        $f_*(N) \le K$.

Then one can prove, in the usual way, that f can be uniquely extended to a commutative cube linking the bicartesian square of s (1.6.1) to the bicartesian square of t:

(2)
$$\begin{array}{ccccc} & & A & \xrightarrow{f} & B \\ & \nearrow^{m} \;\;{\scriptstyle q} & & \nearrow^{h} & \downarrow v \\ M & \xrightarrow{} & & H & \\ q' \downarrow & A/N & \dashrightarrow & B/K & \\ & \nearrow_{m'} & v' \downarrow & \nearrow_{h'} & \\ M/N & \xrightarrow{g} & H/K & & \end{array}$$

This determines the (regularly) *induced* morphism g: M/N → H/K. Notice that, if f is an isomorphism (or even an identity) and the induced morphism g is an isomorphism, its inverse $g^{-1}$ *need not be regularly induced by* $f^{-1}$. This problem is at the basis of the coherence problems dealt with in [G2-G4], in the p-exact setting. Here, we will only deal with regular induction, and induced morphism will mean regularly induced; the general problem of coherence of induction in homological categories (already present in the definition of the derived couple, in 3.2) should be studied.

The commutative cube (2) can be drawn as the following *inductive square* (a cell of the double category Ind**A**, cf. [G7], 2.8)

$$
\begin{array}{ccc}
A & \xrightarrow{f} & B \\
s \uparrow & & \uparrow t \\
M/N & \xrightarrow{g} & H/K
\end{array}
\qquad (3)
$$

Regular induction is consistent with composition and identities. In the normal factorisation of the map $f: A \to B$ (1.2), the central morphism $g$ is induced by $f$ in this sense, from the normal quotient Ncmf (of A), to the normal subobject Nimf (of B).

In particular, if $A = B$ and $f$ is the identity, there are *canonical* morphisms, regularly induced by $1_A$:

(4)    $g: M/N \to M'/N'$,           for $M \leq M'$, $N \leq N'$ in Nsb$A$.

With the previous notation, the direct image of $x \in$ Nsb(M/N) and the inverse image of $y \in$ Nsb(H/K) along the induced morphism $g$ can be computed by direct and inverse images along the edges of the bicartesian squares of $s$ and $t$. This gives four formulas, for each case (as proved in [G7], 1.8):

(5)    $g_*(x) = v'_* h^* f_* m_* q'^*(x) = v'_* h^* f_* q^* m'_*(x) = h'^* v_* f_* m_* q'^*(x) = h'^* v_* f_* q^* m'_*(x)$,

(6)    $g^*(y) = q'_* m^* f^* h_* v'^*(y) = q'_* m^* f^* v^* h'_*(y) = m'^* q_* f^* h_* v'^*(y) = m'^* q_* f^* v^* h'_*(y)$.

It follows that the normal factorisation of the morphism $g$, induced by $f$, can be written as follows

$$
\begin{array}{ccccccc}
(M \wedge f^*K)/N & \xrightarrowtail{k} & M/N & \xrightarrow{g} & H/K & \xtwoheadrightarrow{c} & H/(K \vee f_*M) \\
& & p \downarrow & & \uparrow m & & \\
& & M/(M \wedge f^*K) & \xrightarrow{g'} & (K \vee f_*M)/K & &
\end{array}
\qquad (7)
$$

where *the kernel and cokernel, the normal image and coimage of* $g$ *are canonical morphisms* between subquotients of A and B, while g' is induced by g, hence also by f.

**1.8. The associated perspective category.** Let **A** be an ex2-category, so that the perspective functor Nsb: **A** → **Ltc** is exact (1.5).

We say that **A** is a *perspective* ex2-category if the functor Nsb: **A** → **Ltc** is faithful. This is the case of **Ltc** itself, whose perspective functor is isomorphic to the identity (see 1.4.4). Every ex2-category has an *associated perspective category* Psp**A**, which is the quotient of **A** modulo the congruence which identifies two parallel morphisms $f, g: A \to B$ of **A** whenever Nsb(f) = Nsb(g), or equivalently $f_* = g_*$: Nsb$A \to$ Nsb$B$. It is again an ex2-category, with respect to the 'projected' null morphisms.

It is easy to see that the canonical projection $\mathbf{A} \to \mathrm{Psp}\mathbf{A}$ is an exact functor, which preserves and reflects normal subobjects, normal quotients and exact sequences. Therefore, if $\mathbf{A}$ is homological, or exact, or modular, so is $\mathrm{Psp}\mathbf{A}$.

The p-exact case has already been studied in [G1]. If $\mathbf{A} = K\text{-}\mathbf{Vct}$ is the abelian category of vector spaces on the commutative field $K$, it is easy to see that $f_* = g_*$ if and only if there exists a non-zero scalar $\lambda \in K$ such that $f = \lambda g$. One can thus view the p-exact category $\mathrm{Psp}(K\text{-}\mathbf{Vct})$ as *the category of projective spaces over* $K$; notice that this category has no products, unless $K$ is the two-element field (in which case $K\text{-}\mathbf{Vct}$ is already perspective). We refer to [G1] for a detailed analysis of these facts, and to [CG] for a characterisation of the 'projective categories' derived from the abelian ones.

Therefore, one can think of $\mathrm{Psp}\mathbf{A}$ as a 'projective category' associated to $\mathbf{A}$, also in the general case of ex2-categories. However, we prefer to keep the term 'perspective', to avoid confusion with notions related to 'projective objects'.

The present approach, which makes a crucial use of the transfer functor Nsb, can thus be seen as a 'perspective view' of non-classical homological algebra, as suggested in the Introduction.

## 2. Examples

The category $\mathbf{Set}_2$ of pairs of sets (2.1) is a homological category, with respect to a natural ideal of null morphisms; it is a basic example, as it is a monoidal closed category over which every semiexact category is naturally enriched (2.2).

The category $\mathbf{Top}_2$ of pairs of topological spaces, in the usual sense of algebraic topology, has a similar homological structure, which agrees with the notion of a relative (co)homology theory, in the sense of Eilenberg-Steenrod (2.3). The same holds for the category $\mathbf{Gp}_2$ of pairs of groups, consistently with relative cohomology of groups (2.5). We end with various other homological categories, that are not Puppe-exact (2.6).

**2.1. A basic example.** The category $\mathbf{Set}_2$ of *pairs of sets* has for objects the pairs $(X, X_0)$ where $X_0$ is a subset of the small set $X$; a morphism $f: (X, X_0) \to (Y, Y_0)$ is a mapping $f$ from $X$ to $Y$ such that $f(X_0) \subset Y_0$; the composition is obvious. The morphism $f: (X, X_0) \to (Y, Y_0)$ is assumed to be *null* whenever $f(X) \subset Y_0$; therefore, the null objects are the pairs $(X, X)$ and $f$ is null if and only if it factorises through $(X, X)$ (or, equivalently, through $(Y_0, Y_0)$).

Kernels and cokernels exist and the normal factorisation of $f$ is:

(1)
$$(f^{-1}(Y_0), X_0) \rightarrowtail (X, X_0) \xrightarrow{f} (Y, Y_0) \twoheadrightarrow (Y, Y_0 \cup f(X))$$
$$p \downarrow \qquad \qquad \uparrow m$$
$$(X, f^{-1}(Y_0)) \xrightarrow{g} (Y_0 \cup f(X), Y_0)$$

whence $f$ is exact if and only if it is injective and $f(X) \supset Y_0$. Every normal subobject and every normal quotient of $(X, X_0)$ is determined by a set $A$ with $X_0 \subset A \subset X$:

(2) $(A, X_0) \rightarrowtail (X, X_0) \twoheadrightarrow (X, A)$,

and every short exact sequence in **Set**$_2$ is - up to isomorhism - of this type, determined by a *triple* of sets $X_0 \subset A \subset X$.

**Set**$_2$ is homological: normal monomorphisms and epimorphisms are plainly stable under composition, whereas the existence of subquotients appears from the diagram:

(3)
$$\begin{array}{ccc} (A, X_0) & \rightarrowtail & (X, X_0) \\ \downarrow & & \downarrow \\ (A, B) & \rightarrowtail & (X, B) \end{array}$$

where $(X, X_0)$ is a pair of sets and $X \supset A \supset B \supset X_0$.

Note the following facts, in contrast with the behaviour of p-exact (or abelian) categories: a null morphism $(X, X_0) \to (Y, Y_0)$, between two given objects, need neither exist (take $X \neq \emptyset$, $Y_0 = \emptyset$), nor be unique. A monomorphism (i.e. an injective mapping of pairs) need not be N-mono. A null morphism need not be exact. An exact monomorphism need not be a normal mono: e.g. the normal quotient $(X, X_0) \twoheadrightarrow (X, A)$ is mono, but is not a normal one (for $X_0 \subset A \subset X$, $X_0 \neq A$). The initial object $(\emptyset, \emptyset)$ and the terminal object $(\{*\}, \{*\})$ are distinct; the fact that they are both null is of interest for extending connected sequences of functors.

The category **Set**• of pointed sets $X = (X, 0_X)$ is easily seen to be p-homological, with the following kernels and cokernels

(4) $\quad f^{-1}\{0_Y\} \rightarrowtail X \to Y \twoheadrightarrow Y/f(X)$.

**Set**• is a *perspective* homological category (1.8). The exact functor

(5) $\quad$ P: **Set**$_2 \to$ **Set**•, $\qquad\qquad$ P(X, X_0) = X/X_0,

induces an equivalence of categories Psp**Set**$_2 \to$ **Set**•. (Notice that the 'correct' definition of $X/X_0$ is the pushout of $\{*\} \leftarrow X_0 \to X$; therefore, $X/\emptyset$ is the set $X$ with an extra base-point added.)

**2.2. The canonical enriched structure.** The category of pairs of sets plays for semiexact categories a role similar to that of **Set** for general categories.

Actually, the category **Set**$_2$ of pairs of sets has a natural symmetric monoidal closed structure (somewhat analogous to the smash product of pointed sets, or spaces):

(1) $\quad (X, X_0) \otimes (Y, Y_0) = (X \times Y, (X \times Y_0) \cup (X_0 \times Y))$,

$\qquad$ Hom((X, X_0), (Y, Y_0)) = (**Set**$_2$((X, X_0), (Y, Y_0)), Nul**Set**$_2$((X, X_0), (Y, Y_0))),

where Nul**Set**$_2$((X, X_0), (Y, Y_0)) = **Set**(X, Y_0) is the set of null morphisms, and the adjunction isomorphism is the usual one:

(2) $\quad$ Hom $(X \otimes Z, Y) \to$ Hom $(X,$ Hom $(Z, Y))$, $\qquad f \mapsto (x \mapsto f(x, \text{-}): Z \to Y)$.

Further, **Set**$_2$ has a classifier of normal subobjects:

(3) $\quad$ t: T = $(\{*\}, \{*\}) \to \Omega = (\{0, 1\}, \{1\})$, $\qquad$ t(*) = 1,

$$\begin{CD}
(Y, X_0) @>>> (\{*\}, \{*\}) \\
@VVV @VV{t}V \\
(X, X_0) @>>> (\{0, 1\}, \{1\})
\end{CD}$$

making it a sort of 'monoidal (non cartesian) quasi-topos'; note that the domain of t is the terminal object T and not the identity of the monoidal product, I = ($\{*\}, \emptyset$).

Every semiexact category **A** *with small hom-sets* has an enriched structure over **Set**$_2$, which determines its null morphisms:

(4)   Hom (A, B)  =  (**A**(A, B), Nul**A**(A, B)).

More generally, such an enrichment on a category **C** (with small hom-sets) corresponds to assigning an ideal of **C**, not necessarily closed.

Let us also note that the tensor product is *consistent with the ideal of null morphisms*:

(5)   f $\otimes$ g  is null whenever f or g is.

As an equivalent property, the object (X, X$_0$) $\otimes$ (Y, Y$_0$) is null if one of the factors is.

**2.3. Pairs of spaces.** The category **Top**$_2$ of *pairs of topological spaces* (X, X$_0$), where X$_0$ is a subspace of X, is commonly used as a basis for homological theories. It has a homological structure similar to that of **Set**$_2$: the normal factorisation, short exact sequences and subquotients are still described as above.

In particular, every short exact sequence in **Top**$_2$ is - up to isomorhism - determined by a *triple* of spaces X$_0 \subset$ A $\subset$ X (with respect to the inclusion of subspaces)

(1)   (A, X$_0$) $\rightarrowtail$ (X, X$_0$) $\twoheadrightarrow$ (X, A).

Therefore, a homology theory for pairs of topological spaces, as defined by Eilenberg-Steenrod [ES], gives a 'connected sequence of functors' H$_n$: **Top**$_2$ $\to$ **Ab**.

In analogy to what we have seen at the end of 2.1, the category **Top.** of pointed topological spaces is p-homological, with the 'same' description of kernels and cokernels as in 2.1.4. But now the exact functor

(2)   P: **Top**$_2$ $\to$ **Top.**,                          P(X, X$_0$)  =  X/X$_0$,

(defined again as the pushout of $\{*\} \leftarrow X_0 \to X$) *does not induce an equivalence of categories* Psp**Top**$_2$ $\to$ **Top.**. As well known in Algebraic Topology, the functor P destroys topological information, unless particular hypotheses are assumed on the pairs (X, X$_0$).

**2.4. Groups.** The category **Gp** of groups is semiexact with respect to the ideal of zero morphisms. The normal subobjects can be identified with the invariant subgroups, while every epimorphism is normal; the normal image is the invariant hull of the image; the exact morphisms are those whose image is a normal subgroup. The category of rings has a similar behaviour.

**Gp** is not a homological category, as the normal monomorphisms are not stable under composition. Further, its exactness properties as defined in 1.3, do not agree with the usual definitions: the consecutive morphisms f, g form an exact sequence in that sense if and only if Nim f, the invariant hull of the image, coincides with Ker f, which is weaker than the usual condition Im f = Ker g (note, however, that a *short* exact sequence in the present sense has the usual meaning).

On the other hand, long exact sequences of groups have a well-known habit of getting out of the category, degenerating into pointed sets and actions of groups. Thus, our approach will be to view the exactness properties of groups in larger categories, as the (homological) categories **Gp**$_2$ of pairs of groups (studied below) and **Act**, of actions of groups on pointed sets (Section 5), where the exactness of a sequence of group-homomorphisms resumes its usual meaning. For instance, the long sequences of homotopy 'objects' for a fibration or a pair of pointed spaces can be viewed as exact sequences in **Act** (Section 5); or also in **Gp**$_2$, under suitable restrictions to path connected spaces (4.2).

**2.5. Pairs of groups.** The category **Gp**$_2$ of *pairs of groups* is also analogous to **Set**$_2$ and **Top**$_2$.

An object is a pair $(S, S_0)$ where $S_0$ is a subgroup of the group $S$; a morphism $f: (S, S_0) \to (T, T_0)$ is a group homomorphism from $S$ to $T$ such that $f(S_0) \subset T_0$; it is assumed to be *null* whenever $fS \subset T_0$. The normal factorisation, short exact sequences and subquotients can be described as in **Set**$_2$ (2.1), by replacing subsets and their union with subgroups and their join. The category is homological and, again, a morphism $f: (S, S_0) \to (T, T_0)$ of **Gp**$_2$ is exact if and only if f is injective and $f(S) \supset T_0$. **Gp**$_2$ has a zero object $(0, 0)$, which is only a very particular null object.

The category **Gp** is a retract of **Gp**$_2$, by means of the adjoint (non exact) N-functors K ⊣ I:

(1)   I: **Gp** → **Gp**$_2$,                    $G \mapsto (G, 0)$,

   K: **Gp**$_2$ → **Gp**,                    $(S, S_0) \mapsto S/\bar{S}_0$,

   **Gp**$(S/\bar{S}_0, G) \cong$ **Gp**$_2((S, S_0), (G, 0))$,         KI $\cong$ 1,

where $\bar{S}_0$ is the invariant closure of $S_0$ in S. By adjunction, the functor I is left exact (preserves kernels) while K is right exact (preserves cokernels).

Embedding **Gp** in **Gp**$_2$ via I, a homomorphism of groups $f: G \to H$ has the following normal factorisation *in* **Gp**$_2$:

(2)
$$\begin{array}{ccccc} f^{-1}(0) & \rightarrowtail & G & \xrightarrow{f} & H & \twoheadrightarrow & (H, f(G)) \\ & & p \downarrow & & \uparrow m & & \\ & & (G, f^{-1}(0)) & \xrightarrow{g} & f(G) & & \end{array}$$

Here, Cok f = (H, f(G)) and Nim f = f(G) do not lose information about the image of f, contrarily to what happens in **Gp**; therefore, the exactness of a sequence of group homomorphisms *in* **Gp**$_2$ has its usual meaning. Notice also that the homomorphism f becomes exact in **Gp**$_2$ if and only if it is injective. (This fact, a disadvantage in the study of spectral sequences, can be corrected in the category of normalised groups **Ngp**, where all homomorphisms of groups become exact; cf. 4.7.)

As considered in [G8], Section 3, the *relative* cohomology groups $H^n(S, S_0)$, defined for pairs of groups $(S, S_0)$ (see [M2, Ta, Ri]), form a sequence of functors

(3)   $H^n$: **Gp**$_2$ → **Ab**,                    $(S, S_0) \mapsto H^n(S, S_0)$.

It is actually an 'exact connected sequence of functors': for every triple of groups $S \supset S_1 \supset S_0$, there is a natural exact cohomology sequence.

**2.6. Other examples.** Many other examples of homological categories, non generalised-exact, are considered in the papers [G6, G8]. For instance:

- Banach spaces, Hilbert spaces; topological (or Hausdorff) vector spaces [G6, 4.2, 6.6];
- general 'categories of pairs', extending **Set**$_2$, **Top**$_2$, **Gp**$_2$ [G6, Chapter 4 and 6.7];
- topological coverings [G8, 2.7];
- sets and partial mappings, a category equivalent to **Set.**;
- topological spaces and partial continuous mappings [G6, 4.4, 6.6];
- locally compact $T_2$-spaces and partial proper maps, defined on open subspaces [G6, 5.5, 6.6].

The category **EX**$_4$ of generalised exact $\mathcal{U}$-categories (for some universe $\mathcal{U}$) and exact functors is (unrestricted) homological, with respect to the ideal of exact functors which annihilate every object; the normal subobjects coincide with the thick subcategories, the normal quotients with generalised-exact categories of fractions. The same holds for the categories of p-exact or abelian $\mathcal{U}$-categories [G6, 9.9].

## 3. The exact couple

The exact couple and its associated spectral sequence, introduced by Massey [M1], are here extended to homological categories (following [G7]). This extension will be applied to *homotopy* spectral sequences, whose terms do not live in any abelian context but can be seen as objects of a suitable homological category.

**A** is always a homological category and **BgrA** the homological category of bigraded objects over **A**, with bidegree in $\mathbb{Z}\times\mathbb{Z}$ and morphisms of every bidegree (r, s)

$$f = (f_{hk}): A \to B, \qquad f_{hk}: A_{hk} \to B_{h+r,k+s}.$$

Of course, f is assumed to be a null morphism of **BgrA** when each component is, in **A**. The use of **BgrA** will simplify the exposition, allowing us to defer to the end the issue of bigraduation.

Notice that, even if **A** is pointed, **BgrA** has no zero-object (generally): the object whose components are zero is just *weakly* initial and terminal; thus, if **A** is p-exact, or even abelian, the category **BgrA** is just generalised exact, in the present sense.

**3.1. The exact couple.** An *exact couple* C = (D, E, u, v, ∂) in the homological category **A** is a system of objects and morphisms

(1)
$$\begin{array}{ccc} D & \xrightarrow{u} & D \\ {}_{\partial}\nwarrow & & \swarrow_{v} \\ & E & \end{array}$$

so that

(a) the triangle (1) is exact, i.e. nim u = ker v, nim v = ker ∂ and nim ∂ = ker u.

(b) all the endomorphisms $u^r = u \circ \ldots \circ u: D \to D$ are exact, for $r \geq 1$,

(c)  v  is left modular (1.5) on  Ker u$^r$,  for  r ≥ 1,

(d)  ∂  is right modular (1.5) on  Nim u$^r$,  for  r ≥ 1.

More generally, a *semiexact couple* in  **A**  will only be required to satisfy (a). Note that the endomorphism  d = v∂: E → E  has a null square,  dd = v(∂v)∂,  as  ∂v  is null.

The conditions (b) - (d) are trivially satisfied in a generalised exact category. They can be viewed as a sort of C$^\infty$ condition on our couple; indeed, if *k-exact* couple means that (a) - (d) are satisfied for  1 ≤ r ≤ k,  it is easy to see that the derived couple  C'  is *(k–1)-exact*, following 3.2, 3.3.

**3.2. The derived couple.** The derived couple  C' = (D', E', u', v', ∂')  of the exact couple  C  has the following objects

(1)  D' = Nim u ∈ Nsb(D),

E' = H(E, d) = ∂*(Ker v) / v$_*$(Nim ∂) = ∂*(Nim u) / v$_*$(Ker u).

In order to define its morphisms we also use an isomorphic copy of  D'  *coming from the exactness of the morphism*  u;  it is a normal quotient of  D  (instead of a normal subobject)

(2)  D̄' = Ncm u ∈ Nqt(D),           i: D' ≅ D̄'     (isomorphism induced by  u: D → D).

Now, we have the following morphisms, regularly induced on subquotients (1.7):

(3)  u': D' → D'          induced by  u,  since:  u$_*$(D') = Nim u$^2$ ≤ D',

∂': E' → D'          induced by  ∂,  since:  ∂$_*$(∂*(D')) ≤ D'  and  ∂$_*$(v$_*$(Ker u)) = 0,

v̄': D̄' → E'          induced by  v,  since:  v$_*$(1) = ∂*(0) ≤ Num E',  v$_*$(Ker u) = Den E',

and

(4)  v' = v̄'.i$^{-1}$: D' → D̄' → E'.

Notice that  i$^{-1}$  and  v'  are *not* regularly induced morphisms, generally (1.7).

**3.3. Theorem** ([G7], Thm. 7.3). The derived couple of an exact couple is exact. If the morphisms  v, ∂  are also exact, so are  v'  and  ∂'.                                                                                                                                                    □

**3.4. Iterated derivation.** Therefore, an exact couple  C  has derived couples of any order. Let us begin by *defining* the *r-th derived couple*  C$^r$ = (D$^r$, E$^r$),  for  r ≥ 1,  as

(1)  D$^r$ = Nim u$^{r-1}$ ∈ Nsb(D),                    D$_r$ = Ncm u$^{r-1}$ ∈ Nqt(D),

i: D$_r$ ≅ D$^r$                                    (isomorphism induced by  u$^{r-1}$ : D → D).

E$^r$ = ∂*(Nim u$^{r-1}$) / v$_*$(Ker u$^{r-1}$) = ∂*(D$^r$) / v$_*$(Ker u$^{r-1}$).

The object  E$^r$  is a subquotient of  E,  since  v$_*$(Ker u$^{r-1}$) ≤ v$_*$(1) ≤ ∂*(0) ≤ ∂*(D$^r$). The morphisms of the r-th derived couple are:

(2)  u$^{(r)}$: D$^r$ → D$^r$,      induced by  u      (u$_*$(D$^r$) = D$^{r+1}$ ≤ D$^r$),

∂$^{(r)}$: E$^r$ → D$^r$,      induced by  ∂      (∂$_*$(∂*(D$^r$)) ≤ D$^r$,  ∂$_*$(v$_*$(Ker u$^{r-1}$)) = 0),

v$_{(r)}$: D$_r$ → E$^r$,      induced by  v      (v$_*$(1) ≤ ∂*(0) ≤ Num E$^r$,  v$_*$(Ker u$^{r-1}$) = Den E$^r$),

$$v^{(r)}: (D^r \cong D_r \to E^r), \qquad\qquad v^{(r)} = v_{(r)}.i^{-1}.$$

Now $C = C^1$, $C' = C^2$ and it suffices to show that, assuming the exactness of $C^r$, the couple $C^{r+1}$ is isomorphic to the derived couple $(C^r)'$. Indeed

(3) $\quad (D^r)' = \text{Nim } u^{(r)} = u^{(r)}{}_*(D^r) = u^{(r)}{}_*(\text{Nim } u^{r-1}) = \text{Nim } u^r = D^{r+1}$,

while $(E^r)'$ is the homology of $E^r$ with respect to $d^r = \partial^{(r)}.v^{(r)}: E^r \to E^r$, whence (using also 1.7)

(4) $\quad \text{Ker } d^r = \partial^{(r)*}(\text{Ker } v^{(r)}) = \partial^{(r)*}(\text{Nim } u^{(r)}) = \partial^{(r)*}(D^{r+1}) = (\text{Num } E^r \wedge \partial^*(D^{r+1}))\,/\,\text{Den } E^r =$

$\qquad = \partial^*(D^{r+1})\,/\,\text{Den } E^r,$

$\quad \text{Nim } d^r = v^{(r)}{}_*(\text{Nim } \partial^{(r)}) = v^{(r)}{}_*(\text{Ker } u^{(r)}) = v_{(r)*}(\text{Ker } \bar u^{(r)}) = v_{(r)*}(\text{Ker } u^r\,/\,\text{Den } E^r) =$

$\qquad = (v_*(\text{Ker } u^r) \vee \text{Den } E^r)\,/\,\text{Den } E^r = v_*(\text{Ker } u^r)\,/\,\text{Den } E^r,$

$\quad (E^r)' = H(E^r, d^r) = \text{Ker } d^r\,/\,\text{Nim } d^r = \partial^*(D^{r+1})\,/\,v_*(\text{Ker } u^r) = E^{r+1}.$

**3.5. The bigraded case.** A *bigraded exact couple* $C = (D, E, u, v, \partial)$ *of type 1* in the homological category **A** is an exact couple in the category **BgrA** (recalled above), where the morphisms $u, v$ and $\partial$ have bidegree $(0, 1), (0, 0)$ and $(-1, -1)$, respectively.

In other words, it is a system of morphisms in **A**, indexed on $n, p \in \mathbb{Z}$

(1) $\quad u = u_{np}: D_{n,p-1} \to D_{np},$

$\qquad v = v_{np}: D_{np} \to E_{np},$

$\qquad \partial = \partial_{np}: E_{np} \to D_{n-1,p-1},$

such that:

(a) the following sequences are exact

(2) $\qquad \ldots E_{n+1,p} \xrightarrow{\partial} D_{n,p-1} \xrightarrow{u} D_{np} \xrightarrow{v} E_{np} \xrightarrow{\partial} D_{n-1,p-1} \ldots$

(b) all the morphisms $u^r_{np} = u_{np} \ldots u_{n,p-r+1}: D_{n,p-r} \to D_{np}$ are exact, for $r \geq 1$,

(c) $v_{np}$ is left modular on $\text{Ker}(u^r_{n,p+r}: D_{np} \to D_{n,p+r})$, for $r \geq 1$,

(d) $\partial_{np}$ is right modular on $\text{Nim}(u^r_{np}: D_{n,p-r} \to D_{np})$, for $r \geq 1$.

**3.6. The spectral sequence.** Therefore the r-th derived couple $C^r = (D^r, E^r, u^{(r)}, v^{(r)}, \partial^{(r)})$ of a bigraded exact couple $C$ consists of the following bigraded objects and morphisms ($r \geq 1$):

(1) $\quad D^r_{np} = \text{Nim } (u^{r-1}: D_{n,p-r+1} \to D_{np}) \in \text{Nsb } D_{np} \qquad\qquad (D^1_{np} = D_{np}),$

$\qquad D^{np}_r = \text{Ncm } (u^{r-1}: D_{np} \to D_{n,p+r-1}) \in \text{Nqt } D_{np}\,,$

$\qquad i: D^{np}_r \cong D^r_{n,p+r-1}, \qquad \text{isomorphism induced by } u^{r-1}: D_{np} \to D_{n,p+r-1}\,,$

(2) $\quad E^r_{np} = \partial^*(D^r_{n-1,p-1})\,/\,v_*\,(\text{Ker } u^{r-1}: D_{np} \to D_{n,p+r-1}) \qquad (E^1_{np} = E_{np}),$

(3) $\quad u^{(r)}_{np}: D^r_{n,p-1} \to D^r_{np} \qquad\text{(induced by } u),$

$\qquad \partial^{(r)}_{np}: E^r_{np} \to D^r_{n-1,p-1} \qquad\text{(induced by } \partial),$

$$v_{(r)}^{np}: D_r^{np} \to E_{np}^r \qquad \text{(induced by } v\text{)},$$

$$v_{np}^{(r)} = (D_{n,p+r-1}^r \cong D_r^{np} \to E_{np}^r), \qquad\qquad v^{(r)} = v_{(r)} \cdot i^{-1}.$$

Note that the couple $C^r$ has morphisms $u^{(r)}$, $v^{(r)}$, $\partial^{(r)}$ of bidegrees $(0, 1)$, $(0, 1-r)$, $(-1, -1)$; it can be called a bigraded exact couple *of type* r.

This gives a spectral sequence $(E_{np}^r, d_{np}^r)$, with differential of degree $(-1, -r)$

(4) $\quad d_{np}^r = v^{(r)} \cdot \partial^{(r)} \colon E_{np}^r \to E_{n-1,p-r}^r$.

## 4. Homotopy spectral sequences for path-connected spaces and normalised groups

The homotopy spectral sequence of a tower of fibrations between *path connected* spaces can be seen in the present setting of homological categories as produced by an exact couple living in the p-homological category **Ngp** of *normalised groups* (4.1), a category of fractions of **Gp**$_2$ in which **Gp** embeds, so that all the homomorphisms of groups become exact morphisms of **Ngp**.

This category **Ngp** has a simple description as a quotient of a concrete category (4.5), which is inspired by the category of *homogeneous spaces* introduced by Lavendhomme [La].

The extension to arbitrary topological spaces will be given in Sections 5, 6.

**4.1. The goal.** We will construct a p-homological category **Ngp** = **Q**/$\mathcal{R}$ of *normalised groups*, as a quotient of a homological category **Q** = **Gp**$_2'$ containing **Gp**$_2$.

We will also see (in 4.6) that **Ngp** is a category of fractions of **Gp**$_2$,

(1) $\quad$ **Ngp** $= \Sigma^{-1}$**Gp**$_2$, $\qquad\qquad$ P: **Gp**$_2 \to \Sigma^{-1}$**Gp**$_2$,

obtained by 'excision' of invariant subgroups. In other words, $\Sigma$ is the set of maps

(2) $\quad$ p: $(S, S_0) \to (T, T_0)$, $\qquad\qquad$ $p(S) = T$, $\quad S_0 = p^{-1}(T_0)$,

and these maps coincide up to isomorphism with the canonical projections $(S, S_0) \to (S/N, S_0/N)$, where N is an invariant subgroup of S contained in $S_0$ (take N = Ker p). Notice that $\Sigma$ does not satisfy the two-out-of-three property (and is thus strictly smaller than the set of maps which P makes invertible).

The canonical functor obtained from the embedding I: **Gp** $\to$ **Gp**$_2$ (2.5.1)

(3) $\quad$ J = PI: **Gp** $\to$ **Ngp**, $\qquad\qquad$ G $\mapsto$ (G, 0),

is a left exact and short exact embedding, whose properties are studied in Theorem 4.7. As a crucial fact for our applications to homotopy spectral sequences, all the morphisms f: G $\to$ H of **Gp** become exact morphisms in **Ngp**, since the central morphism of the normal factorisation of If in **Gp**$_2$ (2.5.2) (G, Ker f) $\to$ (f(G), 0) is always in $\Sigma$.

It is also possible to rewrite this section using the calculus of fractions of Gabriel-Zisman ([GZ], I.2), not directly in **Gp**$_2$, where the hypotheses are not satisfied, but in an auxiliary quotient (like, for instance, in the construction of the derived category of an abelian one).

(One considers the congruence $f \mathcal{R} g$ of $\mathbf{Gp_2}$ given, for $f, g: (S, S_0) \to (T, T_0)$, by $fs - gs \in T_0$ for all $s \in S$. This congruence is *implicit* in $\Sigma$, in the sense that every functor defined on $\mathbf{Gp_2}$, which makes all $\Sigma$-maps invertible, identifies all pairs of $\mathcal{R}$-equivalent maps - which can be proved as below, in 4.6.3-4.6.5. Then the category of fractions of $\mathbf{Gp_2}$ with respect to $\Sigma$ trivially coincides with the category of fractions of $\mathbf{Gp_2}/\mathcal{R}$, with respect to the image $\bar\Sigma$ of $\Sigma$; and it is easy to see that the latter category of fractions has a *right* calculus.)

**4.2. The homotopy exact couple of a tower of fibrations.** The motivation for our construction is clear. Let us start from a tower of fibrations of *pathwise connected* pointed spaces ([BK], p. 258)

(1) $\quad \ldots\ X_s \xrightarrow{f_s} X_{s-1} \to \ldots X_0 \xrightarrow{f_0} X_{-1} = \{*\}$

and write $i_s: F_s \to X_s$ the fibre of the fibration $f_s: X_s \to X_{s-1}$. (The fibre is not assumed to be path connected, of course.)

Consider the exact homotopy sequence of $f_{-p}: X_{-p} \to X_{-p-1}$, for $p \leq 0$, in $\mathbf{Gp_2}$

$$
\begin{array}{cccccccccccc}
\ldots \pi_{n+1}X_{-p} & \to & \pi_{n+1}X_{-p-1} & \to & \pi_n F_{-p} & \to & \pi_n X_{-p} & \ldots & \pi_1 X_{-p-1} & \to & (\pi_1 X_{-p-1}, H_p) & \to 0 \\
\| & & \| & & \| & & \| & & \| & & \| & \\
 & u & & v & & \partial & & & & & & \\
\ldots D_{n,p-1} & \to & D_{np} & \to & E_{np} & \to & D_{n-1,p-1} & \ldots & D_{0p} & \to & E_{0p} & \to 0
\end{array}
$$
(2)

where $H_p = (f_{-p})_*(\pi_1 X_{-p})$.

All these sequences produce an *exact couple* of $\mathbf{Ngp}$ (with indices $n \geq 0 \geq p$):

(3) $\quad D_{np} = D^1_{np} = \pi_{n+1}X_{-p-1},$

$\quad E_{np} = E^1_{np} = \pi_n F_{-p}, \qquad E_{0p} = (\pi_1 X_{-p-1}, (f_{-p})_*(\pi_1 X_{-p})) \qquad (n > 0),$

(4) $\quad u_{np} = \pi_{n+1}(f_{-p}): \pi_{n+1}X_{-p} \to \pi_{n+1}X_{-p-1},$

$\quad v_{np}: \pi_{n+1}X_{-p-1} \to \pi_n F_{-p}, \qquad v_{0p}: \pi_1 X_{-p-1} \to (\pi_1 X_{-p-1}, (f_{-p})_*(\pi_1 X_{-p})) \qquad (n > 0),$

$\quad \partial_{np} = \pi_n(i_{-p}): \pi_n F_{-p} \to \pi_n X_{-p}.$

The undefined objects are the zero object of $\mathbf{Ngp}$.

Indeed, all these morphisms are group homomorphisms (embedded in $\mathbf{Ngp}$) or normal epimorphisms (all $v_{0p}$); therefore *all the morphisms are exact in* $\mathbf{Ngp}$ (4.1), and the same is true of all the compositions of morphisms $u_{np}$, consistently with our definition of an exact couple in a homological category (3.1).

Of course it is important to know that $J: \mathbf{Gp} \to \mathbf{Ngp}$ is an embedding, so that we are not losing essential information about homotopy groups.

**4.3. Additive combinations.** We now prepare the ground for a concrete construction of $\mathbf{Ngp}$. Some notation concerning general groups will be useful.

Finite $\mathbb{Z}$-linear combinations $\Sigma_i \lambda_i s_i$ for an abelian group $S$ also make sense for an arbitrary group in additive notation; of course we are no longer allowed to reorder terms, and repetitions of the elements $s_i$ must be permitted. It is thus simpler to consider *additive combinations* $\Sigma_i \varepsilon_i s_i$ where $\varepsilon_i$ stays for $\pm 1$ and the index $i$ varies in a finite *totally ordered* set, say $\{1, \ldots n\}$. Therefore

$$-(\Sigma_i \, \varepsilon_i s_i) \;=\; \Sigma_i \, (-\varepsilon_{n-i}).s_{n-i},$$

and the empty additive combination gives the identity $0$ of the group.

(a) If $S$ is a group, the subgroup $\langle X \rangle$ spanned by a subset $X$ is the set of additive combinations $\Sigma_i \, \varepsilon_i x_i$, with $x_i \in X$.

(b) A mapping $f: S \to T$ between two groups is a homomorphism if and only if it preserves all the additive combinations $(f(\Sigma \, \varepsilon_i s_i) = \Sigma \, \varepsilon_i.fs_i)$, if and only if

(1) $\quad (\Sigma \, \varepsilon_i s_i = 0) \;\Rightarrow\; (\Sigma \, \varepsilon_i.fs_i = 0), \qquad$ for $s_i \in S$, $\varepsilon_i = \pm 1$.

(c) The free group generated by a set $X$ can be described as the set $FX$ of formal additive combinations $\Sigma \, \varepsilon_i \hat{x}_i$ of the elements of $X$, provided that two such formulas are identified when they have the same *reduced* combination, obtained by suppressing all occurrences of type $+\hat{x} - \hat{x}$ or $-\hat{x} + \hat{x}$; the sum in $FX$ is obvious.

(d) If $S$ is a group, we write $ES = F|S|$ the free group generated by the underlying set $|S|$ and

(2) $\quad e = e_S \colon ES \to S, \qquad\qquad e(\Sigma \, \varepsilon_i \hat{s}_i) = \Sigma \, \varepsilon_i s_i,$

the canonical homomorphism, given by the *evaluation* of a formal additive combination as an actual additive combination in $S$. (As well-known, these homomorphisms $e_S \colon F|S| \to S$ form the counit of the adjunction between the free-group functor $F$ and the forgetful functor $|\text{-}|\colon \mathbf{Gp} \to \mathbf{Set}$; the unit is given by the embeddings $i_X \colon X \to |FX|$, $x \mapsto \hat{x}$.)

(e) Every pair of groups $(S, S_0)$ is linked to a *free* pair $(ES, \overline{S}_0)$ by a $\Sigma$-map

(3) $\quad e \colon (ES, \overline{S}_0) \to (S, S_0), \qquad\qquad \overline{S}_0 = e^{-1}(S_0) = \{\Sigma \, \varepsilon_i \hat{x}_i \in ES \mid \Sigma \, \varepsilon_i x_i \in S_0\}.$

**4.4. Quasi-homomorphisms.** Consider the category $\mathbf{Q} = \mathbf{Gp}_2'$ of *pairs of groups* (the objects of $\mathbf{Gp}_2$) and *quasi-homomorphisms* $f \colon (S, S_0) \to (T, T_0)$.

By definition, the latter are *mappings* $f \colon |S| \to |T|$ between the underlying sets, such that the following equivalent conditions hold (same notation as in 4.3, with $s, s', s_i \in S$ and $\varepsilon_i = \pm 1$):

(a) $(\Sigma \, \varepsilon_i s_i \in S_0) \;\Rightarrow\; (\Sigma \, \varepsilon_i.fs_i \in T_0),$

(b) $Ef(\overline{S}_0) \subset \overline{T}_0,$

(c) $fS_0 \subset T_0, \qquad f(\Sigma \, \varepsilon_i s_i) - \Sigma \, \varepsilon_i.fs_i \in T_0,$

(c') $fS_0 \subset T_0, \qquad -f(\Sigma \, \varepsilon_i s_i) + \Sigma \, \varepsilon_i.fs_i \in T_0,$

(d) $fS_0 \subset T_0, \qquad f(s + \varepsilon s') - \varepsilon fs' - fs \in T_0,$

(d') $fS_0 \subset T_0, \qquad -f(\varepsilon s + s') + \varepsilon fs + fs' \in T_0.$

(In fact, (a) $\Leftrightarrow$ (b) is obvious. For (a) $\Rightarrow$ (d), take $\Sigma \, \varepsilon_i s_i = (s + \varepsilon s') - \varepsilon s' - s = 0 \in S_0$. Then (d) $\Rightarrow$ (c) can be proved by induction on the length of additive combinations: if (c) holds for $\Sigma \varepsilon_i s_i$, there exist $t, t' \in T_0$ such that $f(\Sigma \, \varepsilon_i s_i + \varepsilon s) = t + f(\Sigma \, \varepsilon_i s_i) + \varepsilon fs = t + t' + \Sigma \, \varepsilon_i f(s_i) + \varepsilon fs$. Finally, for (c) $\Rightarrow$ (a), let $s = \Sigma \, \varepsilon_i s_i \in S_0$. Then $fs \in T_0$ and $\Sigma \, \varepsilon_i.fs_i = (\Sigma \, \varepsilon_i.fs_i - f(\Sigma \, \varepsilon_i s_i)) + fs \in T_0$. In the same way one proves that (a) $\Rightarrow$ (d') $\Rightarrow$ (c') $\Rightarrow$ (a).)

Furthermore, it is easy to see that $f^{-1}T_0$ is a subgroup of $S$: if all $s_i \in f^{-1}T_0$, then $f(\Sigma \, \varepsilon_i s_i) \in T_0 + \Sigma \, \varepsilon_i.fs_i = T_0$.

**Gp$_2$** is a subcategory of **Q**, and the embedding of **Gp** in **Q** is *full* (every quasi-homomorphism with values in a pair $(T, 0)$ is a homomorphism).

The map $f: (S, S_0) \to (T, T_0)$ of **Q** is assumed to be null if the *set* $fS$ is contained in $T_0$. **Q** is semiexact, with the following normal factorisation

(1)
$$(f^{-1}T_0, S_0) \rightarrowtail (S, S_0) \xrightarrow{f} (T, T_0) \twoheadrightarrow (T, \langle T_0 \cup fX \rangle)$$
$$q \downarrow \qquad \qquad \uparrow n$$
$$(S, f^{-1}T_0) \xrightarrow{g} (\langle T_0 \cup fS \rangle, T_0)$$

where $\langle T_0 \cup fS \rangle$ denotes the subgroup of $T$ spanned by the sub*set* $T_0 \cup fS$.

The normal subobjects and normal quotients of $(S, S_0)$ in **Q** are the same as in **Gp$_2$**, determined by subgroups $M$ with $S_0 \subset M \subset S$; every short exact sequence of **Q** is of the following type, up to isomorphism:

(2) $(M, S_0) \rightarrowtail (S, S_0) \twoheadrightarrow (S, M)$ $\qquad (S_0 \subset M \subset S)$.

Therefore, **Q** is a homological category, with the same description of subquotients as in **Gp$_2$**; the latter is a homological subcategory of **Q**, and contains all the isomorphisms of **Q**.

The properties (d) and (d') show that **Q** is contained in the intersection $\mathcal{P}_g \cap \mathcal{P}_d$ of two (homological) categories considered by Lavendhomme [La], the 'left-extended category of group pairs' $\mathcal{P}_g$ (defined as above by the properties $fS_0 \subset T_0$, $f(s+s') - fs' - fs \in T_0$) and the 'right-extended category of group pairs' $\mathcal{P}_d$ (defined by the properties $fS_0 \subset T_0$, $-f(s + s') + fs + fs' \in T_0$). It can be noted that $\mathcal{P}_g$ and $\mathcal{P}_d$ are isomorphic, by associating to each group the opposite one.

**4.5. Normalised groups.** The category **Ngp** = **Q**/$\mathcal{R}$ of *normalised groups* is the quotient up to the congruence of categories $f \mathcal{R} g$ defined by the equivalent conditions (for $f, g: (S, S_0) \to (T, T_0)$ in **Q**)

(a) for every $s \in S$, $fs - gs \in T_0$,

(b) for every $s \in S$, $-fs + gs \in T_0$,

(c) for all $s_i \in S$ and $\varepsilon_i = \pm 1$, $\Sigma \varepsilon_i . fs_i - \Sigma \varepsilon_i . gs_i \in T_0$.

It suffices to show that (a) $\Rightarrow$ (c); using the quasi-homomorphism property of $f$ and $g$, there are some $t, t' \in T_0$ such that

$$\Sigma \varepsilon_i . fs_i - \Sigma \varepsilon_i . gs_i = t + f(\Sigma \varepsilon s_i) - (t' + g(\Sigma \varepsilon_i s_i)) = t + (f(\Sigma \varepsilon_i s_i) - g(\Sigma \varepsilon_i s_i)) - t' \in T_0.$$

A map $[f]$ in **Ngp** is assumed to be null if and only if $f$ is null in **Q**, independently of the choice of a representative. The null objects are the pairs $(S, S)$, as in **Q** and **Gp$_2$**. But **Ngp** is pointed, with zero object $0 = (0, 0) \cong (S, S)$, since the map $0: (S, S) \to (S, S)$ is $\mathcal{R}$-equivalent to the identity, and $(0, 0)$ is initial and terminal in **Q**.

**Ngp** has kernels and cokernels, which have the same description as in **Q** (independently of the representative we choose for $[f]$).

Therefore **Ngp** is p-homological and the canonical functor

(1) $P: \mathbf{Gp_2} \to \mathbf{Ngp}$, $\qquad (S, S_0) \mapsto (S, S_0)$, $\quad f \mapsto [f]$,

given by the composition $\mathbf{Gp}_2 \to \mathbf{Q} \to \mathbf{Ngp}$ is exact, nsb-faithful and nsb-full (1.5).

(Again, **Ngp** is a subcategory of the Lavendhomme category $\mathcal{H}_g = \mathcal{P}_g/\mathcal{R}_g$ of *left homogeneous spaces* obtained by the congruence $\mathcal{R}_g$ described in (a); it is also a subcategory of the category $\mathcal{H}_d = \mathcal{P}_d/\mathcal{R}_d$ of *right homogeneous spaces*, obtained by the congruence $\mathcal{R}_d$ described in (b). Indeed, both congruences restrict to our $\mathcal{R}$ over **Q**. See [La], Section 3.)

**4.6. Theorem.** The functor P: $\mathbf{Gp}_2 \to \mathbf{Ngp}$ 'is' the category of fractions $\Sigma^{-1}\mathbf{Gp}_2$, i.e. it solves the universal problem of making the morphisms of $\Sigma$ invertible, within (arbitrary) categories and functors.

Furthermore, P is exact and also solves this universal problem within semiexact categories and exact functors (or left exact, or right exact, or short exact functors).

**Proof.** (a) First we prove that P carries each map of $\Sigma$ to an isomorphism of **Ngp**.

Given p: $(S, S_0) \to (S/N, S_0/N)$ in $\Sigma$, choose a *mapping* j: $S/N \to S$ such that $p.j = 1$. Then j: $(S/N, S_0/N) \to (S, S_0)$ is a morphism of **Q**, since it satisfies 4.4(a): if $\Sigma \varepsilon_i.p(s_i) \in S_0/N$, then

$$p(\Sigma \varepsilon_i.jp(s_i)) \in S_0/N, \qquad \Sigma \varepsilon_i.jp(s_i) \in S_0.$$

Further, $jp \mathcal{R} 1$ (and $[j][p] = 1$) because, for every $s \in S$, we have $jp(s) - s \in N \subset S_0$.

(b) Every functor G: $\mathbf{Gp}_2 \to \mathbf{C}$ which makes each map of $\Sigma$ invertible in **C** can be uniquely extended to **Ngp**.

Indeed, given a **Q**-map f: $(S, S_0) \to (T, T_0)$, write, as in 4.3, $ES = F|S|$, e: $ES \to S$ the canonical evaluation epimorphism and $\overline{S}_0 = e^{-1}(S_0) = \{\Sigma \varepsilon_i \hat{x}_i \mid \Sigma \varepsilon_i x_i \in S_0\}$. Then the group homomorphism f': $ES \to T$ defined by the mapping |f| gives a map of $\mathbf{Gp}_2$ and the following commutative diagram in **Q**, with $e \in \Sigma$

(1)
$$\begin{array}{ccc} (S, S_0) & \xrightarrow{f} & (T, T_0) \\ e \uparrow \; \nearrow f' & & \\ (ES, \overline{S}_0) & & \end{array} \qquad f' = e_T.Ef,$$

$$f'(\overline{S}_0) = e_T.Ef(\overline{S}_0) \subset e_T(\overline{T}_0) = e_T e_T^{-1}(T_0) = T_0,$$

$[f] = [f'].[e]^{-1}$, in **Ngp**.

Therefore any functor G' which extends G on **Ngp** is uniquely determined, as follows:

(2)  G': $\mathbf{Ngp} \to \mathbf{C}$, $\qquad G'(S, S_0) = (S, S_0), \qquad G'[f] = (Gf).(Ge)^{-1}.$

Let us prove that G' is well defined by these formulas, and indeed a functor. Firstly, we verify that $f \mathcal{R} g$ in **Q** implies $Gf' = Gg'$. Let

(3)  $U = \langle fS \cup gS \rangle$, $\qquad U_0 = T_0 \cap U,$

$N = \langle \Sigma \varepsilon_i.fs_i - \Sigma \varepsilon_i.gs_i \mid s_i \in S, \varepsilon_i = \pm 1 \rangle.$

Now, N is a subgroup of $U_0$, invariant in U, as follows from the following computation (together with the similar one concerning the inner automorphism produced by $\varepsilon.gs$ instead of $\varepsilon.fs$)

(4)  $\varepsilon.fs + (\Sigma \varepsilon_i.fs_i - \Sigma \varepsilon_i.gs_i) - \varepsilon.fs = \varepsilon.fs + \Sigma \varepsilon_i fs_i - \Sigma \varepsilon_i gs_i - \varepsilon.gs + \varepsilon.gs - \varepsilon.fs =$

$= (\varepsilon.fs + \Sigma \varepsilon_i.fs_i) - (\varepsilon.gs + \Sigma \varepsilon_i.gs_i) + \varepsilon.gs - \varepsilon.fs \in N.$

Considering the following diagram in **Gp**$_2$, where f" and g" are the restrictions of f' and g', the arrow Gs is an isomorphism and sf" = sg"

(5)
$$\begin{array}{ccccc} (ES, \bar{S}_0) & = & (ES, \bar{S}_0) & & \\ f' \downarrow \downarrow g' & & f'' \downarrow \downarrow g'' & & \\ (T, T_0) & \xleftarrow{i} & (U, U_0) & \xrightarrow{s} & (U/N, U_0/N) \end{array}$$

we get that Gf' = Gg'.

Secondly, the following computations show that G' is a functor

(6)
$$\begin{array}{ccccc} (S, S_0) & \xdashrightarrow{[f]} & (T, T_0) & \xdashrightarrow{[g]} & (U, U_0) \\ e \uparrow \quad \nearrow_{f'} & & e \uparrow \quad \nearrow_{g'} & & e \uparrow \\ (ES, \bar{S}_0) & \xrightarrow{Ef} & (ET, \bar{T}_0) & \xrightarrow{Eg} & (EU, \bar{U}_0) \end{array}$$

(gf)' = $e_U$.E(gf) = $e_U$.Eg.Ef,

G'[gf] = (G(gf')).(Ge$_S$)$^{-1}$ = Ge$_U$.GEg.GEf.(Ge$_S$)$^{-1}$ =

= (Ge$_U$.GEg.(Ge$_T$)$^{-1}$).((Ge$_T$).GEf.(Ge$_S$)$^{-1}$) = G'[g].G'[f].

(c) We already know (4.5) that P is exact, nsb-faithful and nsb-full. The last assertion now follows easily from the fact that every map of **Ngp** factorises as φ = Pf.(Pe)$^{-1}$ for some f in **Gp**$_2$ and some e ∈ Σ (cf. (1)). Indeed, if **C** is semiexact and the functor G: **Gp**$_2$ → **C** is left exact, take m = ker f in **Gp**$_2$; then

(7)  Pm = ker Pf,              Ps.Pm = ker Pf.(Ps)$^{-1}$ = ker φ,

G'(ker φ) = Gs.Gm = kerGf.(Gs)$^{-1}$ = kerG'(Pf.(Ps)$^{-1}$) = ker G'(φ).

The right-exact case is proved similarly, as well as the short-exact case (taking into account the fact that P is nsb-full). □

**4.7. Theorem** (Normalising groups). The canonical functor

(1)  J = PI: **Gp** → **Ngp**,          S ↦ (S, 0),   f ↦ [f],

is an embedding. It satisfies the following properties

(a) its codomain is semiexact,

(b) it is short exact,

(c) it takes every monomorphism to a normal mono.

Moreover, J is universal for such properties, i.e. every functor G: **Gp** → **B** which also satisfies (a) - (c) factorises uniquely as G = G'J, where G': **Ngp** → **B** is a short exact functor.

Further, in the presence of (a) and (b), the property (c) is equivalent to:

(c') J is left exact and takes *every* morphism to an exact morphism.

**Proof.** J is trivially an embedding, because of the definition of the congruence $\mathcal{R}$ in 4.5. The property (b) is an easy consequence of Theorem 4.6. Indeed, a short exact sequence of **Gp**, say $N \rightarrowtail G \twoheadrightarrow G/N$, is transformed by I into the following (solid) left exact sequence of **Gp**$_2$

(2)
$$(N, 0) \rightarrowtail (G, 0) \rightarrow (G/N, 0)$$
$$\searrow \quad \uparrow p$$
$$(G, N)$$

The latter becomes short exact in **Ngp**, because the projection p belongs to $\Sigma$.

As to (c), all the monomorphisms of **Gp** become normal monos in **Gp**$_2$, and are preserved as such by the exact functor P: **Gp**$_2 \to$ **Ngp**.

Now, let the functor G: **Gp** $\to$ **B** satisfy (a) - (c). It is easy to see that G extends uniquely to a short exact functor $G^+$: **Gp**$_2 \to$ **B** defined by $G^+(S, S_0) = \mathrm{Cok}_\mathbf{B}(G(S_0 \to S))$. Furthermore, an arbitrary short exact sequence of groups $(m, p) = (N \rightarrowtail S \twoheadrightarrow S/N)$ produces in **Gp**$_2$ a diagram

(3)
$$\begin{array}{ccccc} (N, 0) & \stackrel{m'}{\rightarrowtail} & (S, 0) & \stackrel{p'}{\rightarrow} & (S/N, 0) \\ \| & & \| & & \uparrow u \\ (N, 0) & \rightarrowtail & (S, 0) & \stackrel{q}{\twoheadrightarrow} & (S, N) \end{array}$$

which is transformed by $G^+$ into a commutative diagram whose *upper* row is exact (since $G = G^+I$ is short exact) as well as the *lower* row (since $G^+$ itself is short exact).

Therefore $G^+$ carries all the maps of type u in (3) to isomorphisms; these maps u are particular maps of $\Sigma$ (those with codomain of type $IT = (T, 0)$). Given now an arbitrary $\Sigma$-map p: $(S, S_0) \to (S/N, S_0/N)$, the commutative diagram of **Gp**$_2$, with short exact rows

(4)
$$\begin{array}{ccccc} (S_0, N) & \rightarrowtail & (S, N) & \twoheadrightarrow & (S, S_0) \\ v\downarrow & & \downarrow u & & \downarrow p \\ (S_0/N, 0) & \rightarrowtail & (S/N, 0) & \twoheadrightarrow & (S/N, S_0/N) \end{array}$$

is transformed by $G^+$ into a commutative diagram of **B**, with short exact rows; since $G^+u$ and $G^+v$ are isomorphisms, so is $G^+p$ (which is the result of applying the cokernel functor Cok: $\mathbf{B}^2 \to \mathbf{B}$ to the map $(v, u)$).

Therefore, there is precisely one short exact functor G': $\Sigma^{-1}\mathbf{Gp}_2 \to \mathbf{B}$ such that $G'P = G^+$. Finally $G'J = G'PI = G^+I = G$. On the other hand, if G' is short exact and $G'J = G$, then G'P is short exact and extends G to **Gp**$_2$, so that G'P coincides with $G^+$ and G' is uniquely determined.

Finally, the last assertion. If (c) holds and $k = \ker f$, $p = \mathrm{cok}\, k$ in **Gp**, then $f = ip$ with i mono; therefore, in **Ngp**, Ji is a normal mono and (Jk, Jp) a short exact sequence. Therefore $Jf = Ji.Jp$ is an exact morphism and $\ker Jf = \ker Jp$ is equivalent to Jk. The converse follows from the fact that any exact morphism with a null kernel is a normal mono. □

**4.8. Normalised rings.** In a similar way we can construct the p-homological category **Nrn** of normalised rings, solving a similar problem for the homological category **Rng** of (non-necessarily unital) rings.

The additive combinations of 4.3 are replaced with non-commutative polynomials, with coefficients in $\mathbb{Z}$ and degree $\geq 1$. More precisely, every set $W$ generates a free ring $FW$, whose elements are non-commutative polynomials over $W$, i.e. $\mathbb{Z}$-linear combinations of elements of the free semigroup generated by $W$. For example

(1)  $p(X, Y, Z) = 2.XYXZ - ZX^2Y - XYZX$,

is a (non-null!) non-commutative polynomial over $W = \{X, Y, Z\}$. If $x, y, z$ belong to a ring $R$, $p(x, y, z)$ denotes of course the element $2.xyxz - zx^2y - xyzx$ of $R$. The evaluation map $e: ER = F|R| \to R$ carries each formal non-commutative polynomial to its 'value' in $R$.

We now define, as in 4.4, **Nrn** = **R**/$\mathcal{R}$. Clearly, **R** is the category of *pairs of rings* and *quasi-homomorphisms* $f: (R, R_0) \to (S, S_0)$, i.e. mappings $f: |R| \to |S|$ between the underlying sets such that the following equivalent conditions hold ($p$ denotes an arbitrary non-commutative polynomial and $x_i \in R$):

(a)  $p(x_1,..., x_n) \in R_0 \Rightarrow p(fx_1,..., fx_n) \in S_0$,

(b)  $Ef(\overline{R_0}) \subset \overline{S_0}$,

(c)  $fR_0 \subset S_0$,     $fp(x_1,..., x_n) - p(fx_1,..., fx_n) \in S_0$.

The relation $f \mathcal{R} g$ means again that $fx - gx \in S_0$, for all $x \in R$.

## 5. Actions and homotopy theory

In order to remove the restriction to *path-connected* spaces, this chapter is concerned with the category **Act** [G8] of actions of groups on pointed sets, or actions for short, in which **Set.**, **Gp** and **Gp**$_2$ embed naturally.

We prove that **Act** is homological and show that the homotopy sequences of a pair of spaces or of a fibration can be interpreted as exact sequences in this category.

**5.1. The category of actions.** An *action* is a pair $(X, S)$ where $S$ is a group (always in additive notation) and $X$ is a pointed set (whose base-point is written $0$ or $0_X$) equipped with a right action of $S$ on $X$, written as a sum $x + s$ ($x \in X$, $s \in S$) and satisfying the usual axioms:

(1)  $x + 0_S = x$,     $(x+s) + s' = x + (s+s')$                             ($x \in X$; $s, s' \in S$).

If $x+s = x'$, we say that the operator $s$ *links* $x$ and $x'$. Notice that the base point is *not* assumed to be fixed under $S$. (One can see in 5.6 and 5.7 that such an assumption would make our applications impossible.) We often write

$S_0 = Fix_S(0_X) = \{s \in S \mid 0_X + s = 0_X\}$,

the subgroup of operators which leave the base point fixed.

A morphism of actions $f = (f', f''): (X, S) \to (Y, T)$ consists of a morphism $f': X \to Y$ of pointed sets and a group-homomorphism $f'': S \to T$ consistent with the former:

(2)  $f'0_X = 0_Y$,     $f''(s+s') = f''s + f''s'$,     $f'(x+s) = f'x + f''s$.

We shall often write fx or fs instead of f'x or f"s. The composition is obvious.

A morphism $f = (f', f")$ is assumed to be *null* whenever $f'$ is a zero-morphism in the category of pointed sets: i.e., $fx = 0_Y$ for all $x \in X$. The null objects (those whose identity is null) are the pairs $(0, S)$, where $0$ is a pointed singleton; the ideal of null morphisms is closed (1.1): $f$ is null if and only if it factorises through the null object $(\{0_Y\}, T_0)$, where $T_0 = \text{Fix}_T(0_Y)$. The pair $(0, 0)$ formed by the null group acting on the pointed singleton is a zero-object (and a particular null object).

The kernel of the morphism $f = (f', f"): (X, S) \to (Y, T)$ is the following embedding:

(3)  $(X_1, S_1) \rightarrowtail (X, S)$,

$X_1 = f^{-1}\{0_Y\} = \text{Ker } f'$,

$S_1 = f^{-1}T_0 = \{s \in S \mid 0_Y + fs = 0_Y\} = \{s \in S \mid X_1 + s \subset X_1\} = \{s \in S \mid X_1 + s = X_1\} =$

$= \{s \in S \mid s \text{ links two points of } X_1\}$.

(Note that $S_1$ is determined by $X_1$.) The cokernel of $f$ is the natural projection:

(4)  $(Y, T) \twoheadrightarrow (Y/R, T)$,

where $R$ is the T-congruence of $Y$ generated by identifying all the elements of $fX$.

This proves that **Act** is semiexact; we need further work to show that it is homological.

**5.2. Normal subobjects and quotients.** A normal subobject of the object $(X, S)$ in **Act** can be characterised as the embedding of a pair $(X_1, S_1)$ where:

(a) $S_1$ is a subgroup of $S$ and $X_1$ is a pointed subset of $X$ stable under $S_1$,

(b) if $s \in S$ links two points of $X_1$ then $s \in S_1$ (*normality* condition);

or equivalently:

(a') $X_1$ is a pointed subset of $X$, and if $s \in S$ links two points of $X_1$ then $X_1 + s \subset X_1$,

(b') $S_1 = \{s \in S \mid X_1 + s \subset X_1\} = \{s \in S \mid \exists \, x, x' \in X_1 \text{ such that } x + s = x'\}$,

or also:

(a") $X_1$ is a pointed subset of $X$ and $S_1$ is a subgroup of $S$,

(b") if $x \in X_1$ and $s \in S$, then: $x + s \in X_1 \Leftrightarrow s \in S_1$.

Indeed, it is easy to verify that these conditions are equivalent and that every kernel is of this type. Conversely, given an action $(X_1, S_1)$ satisfying these conditions, consider the following relation in the set $X$:

(1)  $x \, R \, x' \Leftrightarrow (x = x' \text{ or } x = x_1 + s, \; x' = x'_1 + s' \text{ with } x_1, x'_1 \in X_1 \text{ and } s - s' \in S_1)$.

It is an equivalence relation: if we also have $x' = x'_2 + t$ and $x" = x_2 + t'$, with $x_2, x'_2 \in X_1$ and $t - t' \in S_1$, then $s' - t \in S_1$ (because it links two points of $X_1$, namely $x'_1$ and $x'_2 = x' - t = x'_1 + s' - t$), whence $s - t' = (s - s') + (s' - t) + (t - t') \in S_1$ and $x = x_1 + s \, R \, x" = x_2 + t'$. Therefore, $R$ is the S-congruence of $X$ generated by identifying all the points in $X_1$; further, $x \, R \, 0_X$ if and only if $x \in X_1$. It follows immediately that the natural projection:

(2)  $p: (X, S) \to (X/R, S)$,            $p(x) = [x]$,    $p(s) = s$,

is a morphism of **Act**, with kernel $(X_1, S_1)$.

It also follows that the cokernel of $m: (X_1, S_1) \to (X, S)$ is $p: (X, S) \to (X/R, S)$, where R is described in (1): indeed, if $f: (X, S) \to (Y, T)$ annihilates on m, then $f'(X_1) = \{0_Y\}$. Thus, we have also determined the normal quotients of $(X, S)$.

By our characterisation, the ordered set $Nsb(X, S)$ of normal subobjects of $(X, S)$ can be identified with the set of parts $X_1 \subset X$ satisfying the condition (a') above. Since this set is plainly closed under arbitrary intersections, $Nsb(X, S)$ is a complete lattice.

**5.3. Normal factorisations.** We can now describe the normal factorisation of f:

(1)
$$\begin{array}{ccccccc} (X_1, S_1) & \rightarrowtail & (X, S) & \xrightarrow{f} & (Y, T) & \twoheadrightarrow & (Y/R', T) \\ & & q \downarrow & & \uparrow n & & \\ & & (X/R, S) & \xrightarrow[g]{} & (Y_1, T_1) & & \end{array}$$

The kernel $(X_1, S_1)$ is computed in 5.1.3; the normal coimage $(X/R, S) = \text{Cok ker } f$ in 5.2.1:

(2) $X_1 = f^{-1}\{0_Y\} = \text{Ker } f'$,

$S_1 = f^{-1}T_0 = \{s \in S \mid 0_Y + fs = 0_Y\} = \{s \in S \mid X_1 + s \subset X_1\} = \{s \in S \mid X_1 + s = X_1\} =$

$= \{s \in S \mid s \text{ links two points of } X_1\}$.

$x \, R \, x' \Leftrightarrow (x = x' \text{ or } x = x_1 + s, \; x' = x_1' + s' \text{ with } x_1, x_1' \in X_1 \text{ and } s - s' \in S_1)$.

The normal image $\text{Nim}(f) = (Y_1, T_1)$ is the least normal subobject of $(Y, T)$ through which f factorises (1.1) and can be characterised as follows:

(i) $T_1$ is the subgroup of T spanned by the elements $t \in T$ which link the elements of $f(X)$, i.e.: $fx + t = fx'$, for some $x, x'$ in X,

(ii) $Y_1 = fX + T_1$.

Indeed this pair $(Y_1, T_1)$ satisfies the conditions 5.2(a), (b) for a normal subobject of $(Y, T)$: $Y_1$ is stable under $T_1$ and, if $y = fx + t_1 \in Y_1$, $t \in T$ and $y+t \in Y_1$, then $y+t = fx' + t_1'$, whence:

(3) $fx + t_1 + t = fx' + t_1'$, $\qquad t_1 + t - t_1' \in T_1$;

it follows that $t \in T_1$. It is now easy to check that $(Y_1, T_1)$ is the least normal subobject through which f factorises.

Last, we know from 5.1 that $\text{Cok } f = (Y/R', T)$ is determined by the least T-congruence R' of Y which identifies the elements of $f'X$; but $\text{Cok } f = \text{Cok } n$, so that R' can be more concretely presented as the congruence determined by $(Y_1, T_1)$:

(4) $y \, R' \, y' \Leftrightarrow (y = y' \text{ or } y = y_1 + t, \; y' = y_1' + t', \text{ with } y_1, y_1' \text{ in } Y_1, \; t, t' \text{ in } T \text{ and } t - t' \in T_1)$

$\Leftrightarrow (y = y' \text{ or } y = fx + t, \; y' = fx' + t', \text{ with } x, x' \text{ in } X, \; t, t' \text{ in } T \text{ and } t - t' \in T_1)$.

On the other hand, since $n = \ker \text{cok } f$:

(5) $T_1 = \text{Fix}_T(0_{Y/R'})$.

**5.4. Theorem.** **Act** is a homological category.

**Proof.** First, the normal subobjects are stable under composition; given

(1)  $(X_2, S_2) \rightarrowtail (X_1, S_1) \rightarrowtail (X, S),$

if $s \in S$ links the elements $x, x' \in X_2 \subset X_1$ ($x+s = x'$), it follows that $s \in S_1$, and then $s \in S_2$.

Second, let the normal quotients $p$ and $q$ be given, with $Y = X/R$ and $Z = Y/R'$

(2)
$$\begin{array}{ccc} (X, S) & \xrightarrow{p} & (Y, S) \\ h \downarrow & & \downarrow q \\ (X/R'', S) & \xrightarrow{g} & (Z, S) \end{array}$$

and let $h = \mathrm{ncm}\, qp = \mathrm{cok}\,\mathrm{ker}\, qp$. There is a unique morphism $g$ which makes the diagram commutative, and we have to show that it is an isomorphism. Let:

(3)  $(X_1, S_1) = \mathrm{Ker}\, qp = (p^{-1}q^{-1}\{0_Z\}, p^{-1}q^{-1}(\mathrm{Fix}_S(0_Z))) = (p^{-1}q^{-1}\{0_Z\}, \mathrm{Fix}_S(0_Z)).$

The invertibility of $g$ amounts to the injectivity of its component $g'$ (since $g'$ is trivially surjective and $g'' = 1_S$), or also to the following condition: if $x, x' \in X$ and $qp(x) = qp(x')$ then $x\, R''\, x'$. Actually, $px\, R'\, px'$, whence $px = y_1 + s$, $px' = y_1' + s'$ with

$y_1, y_1' \in Y_1 = q^{-1}(\{0_Z\}) = p(X_1),$  $\qquad s - s' \in \{s \in S \mid 0_Z + s = 0_Z\} = S_1.$

Therefore $y_1 = px_1$, $y_1' = px_1'$, and $x\, R\, (x_1+s)\, R''\, (x_1' + s')\, R\, x'$, so that: $x\, R''\, x'$.

As to the last point, namely the existence of subquotients, let us start from the central row of the following diagram:

(4)
$$\begin{array}{ccccc} (X_2, S_2) & \xrightarrow{n} & (X, S) & \xrightarrow{p} & (X/R_1, S) \\ u \downarrow & & \| & & \uparrow v \\ (X_1, S_1) & \xrightarrow{m} & (X, S) & \xrightarrow{q} & (X/R_2, S) \\ h \downarrow & & & & \uparrow k \\ (X_1/\bar{R}_2, S_1) & & \xrightarrow{f} & & (Y, S') \end{array}$$

with $m \geq n = \mathrm{ker}\, q$ and $q \geq p = \mathrm{cok}\, m$. Let $n = mu$ and $p = vq$, so that $u = \mathrm{ker}\, qm$ and $v = \mathrm{cok}\, qm$; let $h = \mathrm{cok}\, u = \mathrm{ncm}\, qm$, $k = \mathrm{ker}\, v = \mathrm{nim}\, qm$ and write $f$ the induced morphism.

We have to prove that $f = (f', f'')$ is an iso. The mapping $f'$ is surjective, as:

$Y = v^{-1}0_{X/R_1} = q(q^{-1}v^{-1}0_{X/R_1}) = qX_1 = \mathrm{Im}\, q'm'.$

It is injective, as: $qm(x_1) = qm(x_1')$ if and only if $x_1\, R_2\, x_1'$, if and only if they are equal or $x_1 = x_2 + s$, $x_1' = x_2' + s'$ with $s - s' \in S_2$, if and only if $x_1\, \bar{R}_2\, x_1'$. Finally, the inclusion $f'': S_1 \to S'$ is the identity, since $S' = v^{-1}\mathrm{Fix}_S(0_{X/R_1}) = \mathrm{Fix}_S(0_{X/R_1}) = S_1$.  $\square$

**5.5. Pointed sets as actions of null groups.** The category **Set.** of pointed sets is equivalent to the full homological subcategory (1.6) **Act$_N$** of **Act** consisting of the actions of null groups (on pointed sets).

Indeed, let us start from the adjoint functors $V \dashv U$

(1) $\quad$ U: **Set.** $\to$ **Act**, $\qquad\qquad$ $U(Z) = (Z, 0)$,

$\qquad$ V: **Act** $\to$ **Set.**, $\qquad\qquad$ $V(X, S) = X/S$.

$\quad$ U associates to a pointed set $Z$ the action on $Z$ of the null group, V associates to an action $(X, S)$ its orbit-set $X/S$, pointed in the orbit $0_X + S$ of the base point, and:

(2) $\quad$ $VU \cong 1$, $\qquad\qquad\qquad\qquad$ **Set.**$(X/S, Z) \cong $ **Act**$((X, S), (Z, 0))$.

$\quad$ Plainly, U restricts to an equivalence **Set.** $\to$ **Act**$_N$ which preserves the null morphisms; therefore, U is exact. V too is exact: given a morphism $g: (X, S) \to (Y, T)$, V preserves its cokernel (by adjointness) and its kernel:

(3) $\quad$ Ker $g = (g^{-1}0_Y, g^{-1}T_0)$, $\qquad$ $V(\text{Ker } g) = (g^{-1}0_Y) / (g^{-1}T_0) = \text{Ker}(Vg: X/S \to Y/T)$.

### 5.6. Pairs of groups as transitive actions.
The full subcategory **Act**$_T$ of **Act** consisting of the *transitive* actions is again a homological subcategory, which we now prove to be equivalent to **Gp**$_2$.

$\quad$ There is an adjoint retraction $F \dashv G$, where F is exact and G is left exact, short exact

(1) $\quad$ F: **Gp**$_2 \to$ **Act**, $\qquad\qquad$ $F(S, S_0) = (|S|/S_0, S)$,

$\qquad$ G: **Act** $\to$ **Gp**$_2$, $\qquad\qquad$ $G(X, S) = (S, S_0)$ $\qquad\qquad\qquad\qquad$ $(S_0 = \text{Fix}_S(0_X))$.

$\quad$ F associates to a pair $(S, S_0)$ the canonical (right) action of the group $S$ over the pointed set $|S|/S_0$ of the right cosets of $S_0$ in $S$: $(S_0 + s) + s' = S_0 + (s+s')$, so that:

(2) $\quad$ $GF(S, S_0) = (S, S_0)$, $\qquad\qquad$ **Act**$((|S|/S_0, S), (Y, T)) \cong$ **Gp**$_2((S, S_0), (T, T_0))$.

$\quad$ Further, **Gp**$_2$ is equivalent to the full homological subcategory **Act**$_T$ of **Act**, by means of the restriction F': **Gp**$_2 \to$ **Act**$_T$ of F. Indeed $F(S, S_0) = (|S|/S_0, S)$ is a transitive action, F' is fully faithful and for every transitive action $(X, S)$, the counit $\varepsilon$ is an isomorphism. This also proves that F is exact (sparing a longer direct proof).

$\quad$ The functor G is not exact, but is left exact (as a right adjoint) and also short exact (i.e. it preserves short exact sequences): given a short exact sequence and its G-image

(3) $\quad$ $(X_1, S_1) \rightarrowtail (X, S) \twoheadrightarrow (X/R, S)$,

$\qquad$ $(S_1, S_0) \rightarrowtail (S, S_0) \twoheadrightarrow (S, \text{Fix}(0_{X/R}))$,

by 5.3.5, we have: $\text{Fix}(0_{X/R}) = S_1$, whence $(S, \text{Fix}(0_{X/R})) = (S, S_1) = \text{Cok}((S_1, S_0) \rightarrowtail (S, S_0))$.

$\quad$ Last, the category **Gp** is a retract of **Act**, as it follows by composing the retractions already considered in 2.5 and above

(4) $\quad$ FI: **Gp** $\to$ **Act**, $\qquad$ $S \mapsto (|S|, S)$ $\qquad\qquad\qquad\qquad\qquad\qquad$ (left exact),

$\qquad$ KG: **Act** $\to$ **Gp**, $\qquad$ $(X, S) \mapsto S/\overline{S}_0$ $\qquad\qquad\qquad\qquad$ (preserves normal quotients),

where FI associates to a group $S$ the canonical (right) action $x+s$ of $S$ over the underlying pointed set $|S|$, while KG associates to an action $(X, S)$ the quotient of $S$ modulo the invariant closure of $S_0 = \text{Fix}_S(0_X)$ in $S$.

$\quad$ Note that, again (as in 2.5), a sequence of groups, embedded via FI, is exact in **Act** if and only if it is exact in the usual sense (im $f = $ ker $g$).

**5.7. Exact homotopy sequences in Act.** If $p: X \to B$ is a Serre fibration of pointed spaces, with fibre $F = p^{-1}\{0_B\}$, its homotopy sequence can be written as a sequence in **Act**

(1)    ... $\pi_1 F \to \pi_1 X \to \pi_1 B \to (\pi_0 F, \pi_1 B) \to \pi_0 X \to \pi_0 B$,

where

- $(\pi_0 F, \pi_1 B)$ is the usual action of the group $\pi_1 B$ over the pointed set $\pi_0 F$,
- all the terms at its left are groups, embedded in **Act** (5.6.4),
- the last two terms are pointed sets, embedded in **Act** (5.5.1).

This sequence is *exact* in **Act**, as follows from the theorem below (5.8).

Similarly, it is easy to verify that the homotopy sequence of a pair $(X, A)$ of pointed spaces can be read as an exact sequence in **Act**:

(2)    ... $\to \pi_1 A \to \pi_1 X \to (\pi_1(X, A), \pi_1 X) \to \pi_0 A \to \pi_0 X \to \pi_0(X, A) \to 0$.

All the terms down to $\pi_1 X$ are groups (embedded in **Act**, by 5.6.4); $\pi_1(X, A)$ is the pointed set of paths $\sigma: I \to X$ with $\sigma(0) \in A$, $\sigma(1) = 0_X = 0_A$, modulo homotopy with first end in $A$ and second end in $0_X$, and the group $\pi_1 X$ acts on the right on $\pi_1(X, A)$ in the natural, standard way. The four last terms are pointed sets (embedded in **Act**, by 5.5.1); in particular, $\pi_0(X, A)$ is the cokernel of the pointed mapping $\pi_0 A \to \pi_0 X$, i.e. the set of path-components of $X$ modulo the relation identifying all the components which intersect $A$ (pointed in this class), while the last object $0$ is the singleton (the zero object of **Set.**).

More generally, given a triple of pointed spaces $B \subset A \subset X$, we get an exact sequence in **Act**:

(3)    ... $\to \pi_2(X, A) \to (\pi_1(A, B), \pi_1 A) \to (\pi_1(X, B), \pi_1 X) \to (\pi_1(X, A), \pi_1 X) \to$

   $\to \pi_0(A, B) \to \pi_0(X, B) \to \pi_0(X, A) \to 0$,

whose terms are groups in degree $\geq 2$, 'general' actions in degree 1 and pointed sets in degree 0 (including the last term, 0).

It can be noted that, if all the spaces are pathwise connected, the last three terms in (2) annihilate and all actions become transitive: by 5.6, the sequence (2) can be realised as an exact sequence in **Gp**$_2$.

**5.8. Theorem** (Exactness from groups to pointed sets). Consider the following sequence in **Act**

(1)    $H \xrightarrow{u} G \xrightarrow{v} S \xrightarrow{f} (X, S) \xrightarrow{g} Y \xrightarrow{h} Z$

where

- H, G, S, u, v are in **Gp**, viewed in **Act** as $(|H|, H)$ etc., by the left exact embedding FI (5.6.4),
- $(X, S)$ is an action and $f = (f', \mathrm{id}S)$ with $f'(s) = 0_X + s$ (for $s \in |S|$),
- Y, Z, h are in **Set.**, viewed in **Act** as $(Y, 0)$ etc., by the exact embedding U (5.5.1)
- $g: X \to Y$ is a map of pointed sets such that $g(x + s) = g(x)$ for all $x \in X$, $s \in S$.

Then:

(a) the sequence is exact in G if and only if Im u = Ker v,

(b) the sequence is exact in  S  if and only if  Im v = Ker f' = $\text{Fix}_S\{0_X\}$,

(c) the sequence is exact in  (X, S)  if and only if  $0_X + S = g^{-1}\{0_Y\}$,

(d) the sequence is exact in  Y  if and only if  $g'(X) = h^{-1}\{0_Z\}$,

(e) the morphism  f  is necessarily exact,

(f) the morphism  g  is necessarily right modular (1.5).

*Note*. The classical properties of exact homotopy sequences coincide with the present ones, except in (X, S),  where the classical property is stronger: for  x, x' ∈ X,  g(x) = g(x')  if and *only if*  x = x' + s  for some  s ∈ S.

**Proof.** Point (a) follows from the last remark in 5.6.

Points (b) and (e) follow from the normal factorisation of  f,  computed as follows (by 5.3)

(2)
$$
\begin{array}{ccccc}
(|S_0|, S_0) \rightarrowtail (|S|, S) & \xrightarrow{f} & (X, S) & \twoheadrightarrow & (X/X_0, S) \\
\downarrow & \uparrow & & & S_0 = \text{Fix}_S\{0_X\}, \\
(|S|/S_0, S) & \rightarrow & (X_0, S) & & X_0 = 0_X + S.
\end{array}
$$

Here,  $X/X_0$  is the pointed set obtained by collapsing the orbit  $X_0 = 0_X + S$  to a (base) point. Moreover,  $|S|/S_0$  denotes the set of right cosets  $S_0 + x$.  Indeed, the relation  R  in  |S|  described in 5.2.1 for the normal coimage of  f  becomes now

(3)   x R x'  ⇔  (x = $x_1$ + s,  x' = $x_1'$ + s'  with  $x_1, x_1' \in S_0$  and  s – s' ∈ $S_0$),

which is easily seen to be equivalent to  x – x' ∈ $S_0$.  (If this is the case, just let  x = 0 + x,  x' = 0 + x'; conversely, given (3), we have  x – x' = $x_1$ + s – s' – $x_1'$ ∈ $S_0$.)

Points (c) and (d) follow from the normal factorisation of  g,  which is:

(4)
$$
\begin{array}{ccccc}
(X_1, S) \rightarrowtail (X, S) & \xrightarrow{g} & (Y, 0) & \twoheadrightarrow & (Y/g'(X), 0) \\
\downarrow & \uparrow & & & X_1 = g^{-1}\{0_Y\} = \text{Ker } g'. \\
(X/X_1, S) & \rightarrow & (g'(X), 0) & &
\end{array}
$$

Finally, for (f), let  Y'  be a pointed subset of  Y  (i.e., a normal subobject in **Set.** and **Act**). Then  $g^{-1}(Y')$  is stable under the action of  S,  and we have:

(5)   $g_* g^*(Y') = g_*(g^{-1}(Y'), S) = g(g^{-1}(Y')) = Y' \cap g(X)$.   □

## 6. Homotopy spectral sequences and normalised actions

We introduce a homological category of fractions  **Nac** = $\Sigma^{-1}$**Act**,  which will be called *the category of normalised actions*, and proved to be adequate for homotopy sequences of (possibly) non path-connected spaces.

**6.1. Quasi-exact couples.** Extending 3.5, a *bigraded quasi-exact couple* $C = (D, E, u, v, \partial)$ *of type 1* in the homological category **A** will be a system of objects and morphisms

(1) $\quad D_{n,p} \quad\quad (n \geq 0,\ p \leq 0), \quad\quad\quad E_{n,p} \quad (n \geq 1,\ p \leq 0),$

$\quad\quad u = u_{np}\colon D_{n,p-1} \to D_{np}, \quad\quad\quad\quad\quad\quad\quad\quad\quad\quad\quad (n \geq 0,\ p \leq 0),$

$\quad\quad v = v_{np}\colon D_{np} \to E_{np}, \quad\quad\quad\quad\quad\quad\quad\quad\quad\quad\quad\quad\ (n \geq 1,\ p \leq 0),$

$\quad\quad \partial = \partial_{np}\colon E_{np} \to D_{n-1,p-1} \quad\quad\quad\quad\quad\quad\quad\quad\quad\quad\ \ (n \geq 1,\ p \leq 0),$

such that:

(a) the following sequences are exact

(2) $\quad \ldots\ E_{n+1,p} \xrightarrow{\partial} D_{n,p-1} \xrightarrow{u} D_{np} \xrightarrow{v} E_{np}\ \ldots\ E_{1p} \xrightarrow{\partial} D_{0,p-1} \xrightarrow{u} D_{0p}$

(b) all the morphisms $u_{np}^r = u_{np} \ldots u_{n,p-r+1}\colon D_{n,p-r} \to D_{np}$ are exact, for $r \geq 1$ and $n > 0$,

(c) $v_{np}$ is left modular on $\mathrm{Ker}(u^r_{n,p+r}\colon D_{np} \to D_{n,p+r})$, for $r \geq 1$,

(d) $\partial_{np}$ is right modular on $\mathrm{Nim}(u^r_{np}\colon D_{n,p-r} \to D_{np})$, for $r \geq 1$.

Apart from the restriction on the indices n, p, the real interest of the extension is that *we are not requiring the exactness of* the last morphism in sequence (2), namely $u_{0p}\colon D_{0,p-1} \to D_{0p}$, because this is not satisfied in our application below (6.2). However, the present hypotheses are sufficient to obtain the associated derived couples and spectral sequence, as in 3.6; indeed, in the construction of the derived couple, we only need the morphism

(3) $\quad v_{np}^{(r)} = (D^r_{n,p+r-1} \cong D^{np}_r \to E^r_{np}),$

for $n \geq 1$. Therefore, we only need the isomorphism $i\colon D^{np}_r \to D^r_{n,p+r-1}$ *induced by the exact morphism* $u^{r-1}\colon D_{np} \to D_{n,p+r-1}$ for $n \geq 1$.

**6.2. A category of fractions.** The category of fractions $\mathbf{Nac} = \Sigma^{-1}\mathbf{Act}$ in which we are interested is obtained by 'excision of the invariant subgroups of operators which act trivially', i.e. the set $\Sigma$ is formed – up to isomorphism – by all the natural projections $p\colon (X, S) \to (X, S/N)$, where $N$ is an invariant subgroup of $S$ which acts trivially on $X$.

The motivation comes from a tower of fibrations of *arbitrary* pointed spaces ([BK], p. 258)

(1) $\quad \ldots\ X_s \xrightarrow{f_s} X_{s-1} \to \ldots X_0 \xrightarrow{f_0} X_{-1} = \{*\}$

Again, we write $i_s\colon F_s \to X_s$ the fibre of the fibration $f_s\colon X_s \to X_{s-1}$.

Consider the exact homotopy sequence of $f_{-p}\colon X_{-p} \to X_{-p-1}$, for $p \leq 0$, in **Act** (5.7)

(2) $\quad \ldots\ \pi_1 F_{-p} \to \pi_1 X_{-p} \to \pi_1 X_{-p-1} \to (\pi_0 F_{-p}, \pi_1 X_{-p-1}) \to \pi_0 X_{-p} \to \pi_0 X_{-p-1}$

$\quad\quad\quad\ \ \|\quad\quad\quad\ \|\quad\quad\quad\ \|\quad\quad\quad\quad\quad\ \|\quad\quad\quad\quad\ \|\quad\quad\quad\ \|$

$\quad\quad\quad\quad\ \ \xrightarrow{\partial}\quad\quad\ \xrightarrow{u}\quad\quad\ \xrightarrow{v}\quad\quad\quad\quad\quad\quad\ \xrightarrow{\partial}\quad\quad\quad\ \xrightarrow{u}$

$\quad \ldots\ E_{2p} \to D_{1,p-1} \to D_{1p} \to E_{1p} \to D_{0,p-1} \to D_{0p}$

where $(\pi_0 F_{-p}, \pi_1 X_{-p-1})$ is the well-known canonical action, the objects at its left are groups (embedded in **Act**, by 5.6.4), and the last two terms are pointed sets (embedded in **Act**, by 5.5.1).

All these sequences produce a semiexact couple in **Act**:

(3) $\quad D_{np} = D^1_{np} = \pi_n X_{-p-1}$, $\hfill (n \geq 0)$,

$\quad E_{np} = E^1_{np} = \pi_{n-1} F_{-p}, \qquad\qquad E_{1p} = (\pi_0 F_{-p}, \pi_1 B) \hfill (n > 1)$,

(4) $\quad u_{np} = \pi_n(f_{-p}): \pi_n X_{-p} \to \pi_n X_{-p-1}, \hfill (n \geq 0)$,

$\quad v_{np}: \pi_n X_{-p-1} \to \pi_{n-1} F_{-p}, \qquad v_{1p}: \pi_1 X_{-p-1} \to (\pi_0 F_{-p}, \pi_1 X_{-p-1}) \hfill (n > 1)$,

$\quad \partial_{np} = \pi_{n-1}(i_{-p}): \pi_{n-1} F_{-p} \to \pi_{n-1} X_{-p}, \quad \partial_{1p} = \pi_0(i_{-p}): (\pi_0 F_{-p}, \pi_1 X_{-p-1}) \to \pi_0 X_{-p}, \hfill (n \geq 2)$.

It becomes a *quasi-exact couple* in **Nac**. Indeed, all these morphisms fall in the following situations:

- group homomorphisms (embedded in **Act** and) projected to exact morphisms in **Nac**,

- morphisms $v_{1p}: \pi_1 X_{-p-1} \to (\pi_0 F_{-p}, \pi_1 X_{-p-1})$ which are exact in (**Act** and) **Nac**, by 5.8(e), whence also left modular,

- morphisms $\partial_{1p} = \pi_0(i_{-p}): (\pi_0 F_{-p}, \pi_1 X_{-p-1}) \to \pi_0 X_{-p}$, which are right modular by 5.8(f).

On the other hand, the morphisms of pointed sets: $u_{0p} = \pi_0(f_{-p}): \pi_0 X_{-p} \to \pi_0 X_{-p-1}$ need not be exact, and we only have a *quasi* exact couple.

The sequel is devoted to construct **Nac** and verify the desired properties.

**6.3. A larger category of actions.** Again, we will construct $\Sigma^{-1}$**Act** as a quotient of a homological category **Act'**, following the same line as for $\mathbf{Q} = \mathbf{Gp}'_2$, in Section 4.

An object of **Act'** is a triple $(X, S, S_0)$, where $(X, S)$ belongs to **Act**, $(S, S_0)$ is a pair of groups, and $S_0$ acts trivially on $X$. A morphism $f = (f', f''): (X, S, S_0) \to (Y, T, T_0)$ consists of a map of pointed sets $f': X \to Y$ and a **Q**-morphism $f'': (S, S_0) \to (T, T_0)$ which are consistent: $f'(x+s) = f'x + f''s$, for $x \in X$ and $s \in S$. Composition is obvious.

The map $f = (f', f''): (X, S, S_0) \to (Y, T, T_0)$ of **Act'** is assumed to be null if $f': X \to Y$ is a zero-map of pointed sets, i.e. $f'(X) = \{0_Y\}$. **Act'** is semiexact, with the following normal factorisation of $f$ (obtained from the normal factorisation of a morphism in **Act**, see 5.3.1)

(1)
$$\begin{array}{ccccccc}
(X_1, S_1, S_0) & \stackrel{m}{\rightarrowtail} & (X, S, S_0) & \stackrel{f}{\longrightarrow} & (Y, T, T_0) & \stackrel{p}{\twoheadrightarrow} & (Y/R', T, T_0) \\
& & q \downarrow & & \uparrow n & & \\
& & (X/R, S, S_0) & \stackrel{}{\underset{g}{\longrightarrow}} & (Y_1, T_1, T_0) & &
\end{array}$$

The definition of $X_1, S_1, R, R', T_1, Y_1$ is the same as in 5.3 (even if $f''$ is just a quasi-homomorphism); notice that $S_1 = \{s \in S \mid X_1 + s \subset X_1\}$ contains $S_0$ and $T_1$ (the subgroup of $T$ spanned by the elements which link the elements of $f'(X)$) contains $T_0$.

One proves that **Act'** is homological, as in 5.4. There is an exact embedding

(2) $\quad$ **Act** $\to$ **Act'**, $\qquad\qquad (X, S) \mapsto (X, S, S_0)$,

where $S_0 = \text{Fix}_S(0_X)$. (Notice that replacing $S_0$ with $\text{Fix}_S(X)$ would not give a functor: from $s \in \text{Fix}_S(X)$ we can only deduce that $f(s) \in \text{Fix}_T(fX)$.)

**6.4. Normalised actions.** The category **Nac** = **Act'**/$\mathcal{R}$ of *normalised actions* is the quotient modulo the congruence of categories $f \mathcal{R} g$ defined by $f' = g'$ and $f'' \mathcal{R} g''$ in **Q** (4.5). The latter amounts to the equivalent conditions

(a) for every $s \in S$, $f''s - g''s \in T_0$,

(b) for every $s \in S$, $-f''s + g''s \in T_0$,

(c) for all $s_i \in S$ and $\varepsilon_i = \pm 1$, $\Sigma \varepsilon_i.f''s_i - \Sigma \varepsilon_i.g''s_i \in T_0$.

A map [f] in **Nac** is assumed to be *null* if and only if f' is null in **Set.**. The null objects are the triples $(0, S, S_0)$.

**Nac** has kernels and cokernels, with the same description as in **Act'** (independently of the representative we choose for [f]).

Therefore **Nac** is homological and the canonical functor

(1)   P: **Act** $\to$ **Nac**,

   $(X, S) \mapsto (X, S, S_0)$,   $f \mapsto [f]$                                                    $(S_0 = \text{Fix}_S(0_X))$,

given by the composition **Act** $\to$ **Act'** $\to$ **Nac** is exact, nsb-faithful and nsb-full (1.5).

**6.5. Theorem.** The functor P: **Act** $\to$ **Nac** 'is' the category of fractions $\Sigma^{-1}$**Act**, i.e., it solves the universal problem of making the morphisms of $\Sigma$ (6.2) invertible (within arbitrary categories and functors).

Furthermore, P is exact and also solves this universal problem within semiexact categories and exact functors (or left exact, or right exact, or short exact functors).

**Proof.** The proof follows the same line as for Theorem 4.6, with suitable modifications.

(a) P carries each map of $\Sigma$ to an isomorphism of **Nac**. Given p: $(X, S) \to (X, S/N)$ in $\Sigma$ (with N a normal subobject of S which acts trivially on X), we want to prove that the associated map

   $\hat{p}$: $(X, S, S_0) \to (X, S/N, S_0/N) \in \text{Mor}\mathbf{Act'}$,                         $(S_0 = \text{Fix}_S(0_X))$,

becomes an isomorphism in **Nac**.

Choose a *mapping* j: S/N $\to$ S such that p.j = 1. Then

   $\hat{j} = (\text{id}X, j)$: $(X, S/N, S_0/N) \to (X, S, S_0)$,

is a morphism of **Act'**, since:

- if $\Sigma \varepsilon_i.p(s_i) \in S_0/N$, then $p(\Sigma \varepsilon_i.jp(s_i)) \in S_0/N$ and $\Sigma \varepsilon_i.jp(s_i) \in S_0$,

- if $x \in X$ and $s \in S$, $x + jp(s) = x + s = x + p(s)$, because $-s + pj(s) \in N \subset \text{Fix}_S(X)$.

By the same reason, $\hat{j}\hat{p} \mathcal{R} 1$.

(b) Every functor G: **Act** $\to$ **C** which makes each $\Sigma$-map invertible in **C** can be uniquely extended to **Nac**. Given an **Act'**-map f = (f', f''): $(X, S, S_0) \to (Y, T, T_0)$ write, as in 4.3, ES = F|S| and e:

$ES \to S$ the canonical evaluation epimorphism; $ES$ acts on $X$ in the obvious way, by evaluating $ES$ in $S$, and $\bar{S}_0 = e^{-1}(S_0) = \{\Sigma\, \varepsilon_i \hat{s}_i \mid \Sigma\, \varepsilon_i s_i \in S_0\}$ acts trivially on $X$.

Then the *group homomorphism* $f_1'' = e_T.Ef\colon ES \to T$ defined by the *mapping* $|f''|\colon |S| \to |T|$ gives a map $f_1 = (f', f_1'')\colon (X, ES) \to (Y, T)$ of **Act**, that is also in **Act**'

(1)
$$
\begin{array}{c}
(X, S, S_0) \xdashrightarrow{f} (Y, T, T_0) \\
e \uparrow \quad \nearrow f_1 \\
(X, ES, \bar{S}_0)
\end{array}
\qquad
\begin{array}{l}
f_1 = (f', f_1'') = (f', e_T.Ef), \\
[f] = [f_1].[e]^{-1}, \text{ in } \mathbf{Nac},
\end{array}
$$

because $f_1''(\bar{S}_0) = e_T.Ef(\bar{S}_0) \subset e_T(\bar{T}_0) = e_T\, e_T^{-1}(T_0) = T_0$,

Therefore any functor $G'$ which extends $G$ on **Nac** is uniquely determined, as follows:

(2)   $G'\colon \mathbf{Nac} \to \mathbf{C}$,     $G'(X, S, S_0) = (X, S, S_0)$,    $G'[f] = (Gf_1).(Ge)^{-1}$.

Let us prove that $G'$ is well defined by these formula, and indeed a functor. Firstly, we verify that $f\, \mathcal{R}\, g$ in **Act**' implies $Gf_1 = Gg_1$. Let

(3)   $T' = \langle fS \cup gS \rangle$,       $Y' = (fX \cup gX) + T'$,

   $T'' = T_0 \cap T'$,       $N = \langle \Sigma\, \varepsilon_i.fs_i - \Sigma\, \varepsilon_i.gs_i \mid s_i \in S,\ \varepsilon_i = \pm 1 \rangle$.

Now, $N$ is a subgroup of $T''$, invariant in $T'$ (with the same proof as in 4.6.4). Therefore the following diagram in **Act**', where $f_2$ and $g_2$ are the restrictions of $f_1$ and $g_2$

(4)
$$
\begin{array}{ccc}
(X, ES, \bar{S}_0) & = & (X, ES, \bar{S}_0) \\
f_1 \downarrow\downarrow g_1 & & f_2 \downarrow\downarrow g_2 \\
(Y, T, T_0) \xleftarrow{i} (Y', T', T'') & \xrightarrow{p} & (Y'/N, T''/N)
\end{array}
$$

shows that $Gf_1 = Gg_1$ (because $Gp$ is an isomorphism and $pf_2 = pg_2$).

Finally, $G'$ is a functor, with the same computations as in 4.6.6, based on a slightly different diagram

(5)
$$
\begin{array}{ccccc}
(X, S, S_0) & \xdashrightarrow{[f]} & (Y, T, T_0) & \xdashrightarrow{[g]} & (Z, U, U_0) \\
e \uparrow \quad \nearrow f_1 & & e \uparrow \quad \nearrow g_1 & & e \uparrow \\
(X, ES, \bar{S}_0) & \xrightarrow[Ef]{} & (Y, ET, \bar{T}_0) & \xrightarrow[Eg]{} & (Z, EU, \bar{U}_0)
\end{array}
$$

(c) Same argument as for point (c) in Theorem 4.6.    □